\documentclass[11pt]
{amsart}
\usepackage{amssymb,amsmath,amsthm,amsfonts,amsopn,url,xcolor,hyperref,enumerate,mathtools,microtype,MnSymbol,comment}
\usepackage[normalem]{ulem}
\usepackage[all]{xy}
\input xy
\xyoption{all}
\usepackage{amscd}
\usepackage{soul}
\usepackage[mathscr]{euscript}

\theoremstyle{plain}
\newtheorem{thm}{Theorem}[section]

\newtheorem{prop}[thm]{Proposition}
\newtheorem{lemma}[thm]{Lemma}
\newtheorem{cor}[thm]{Corollary}

\theoremstyle{definition}
\newtheorem{defn}[thm]{Definition}
\newtheorem*{defn*}{Definition}
\newtheorem*{question*}{Question}

\newtheorem{example}[thm]{Example}
\newtheorem*{example*}{Example}
\newtheorem{rem}[thm]{Remark}
\newtheorem*{rem*}{Remark}

\newtheorem{notation}[thm]{Notation}
\newtheorem{algorithm}[thm]{Algorithm}
\newcommand{\field}[1]{\mathbb{#1}}

\newcommand{\R}{\field{R}}
\newcommand{\F}{\field{F}}

\newcommand{\func}[1]{\mathrm{#1} \,}

\newcommand{\im}{\func{im}}

\newcommand{\be}{\begin{enumerate}}
\newcommand{\ee}{\end{enumerate}}

\newcommand{\li}
 {\leftfootline}

\newcommand{\cA}{\mathcal{A}}

\newcommand{\cP}{\mathcal{P}}
\renewcommand{\phi}{\varphi}

\let\int\relax
\DeclareMathOperator{\int}{i}

\newcommand{\fourstep}{recomposed}
\newcommand{\fourstepping}{recomposing}
\newcommand{\fivestep}{almost canonical form}

\author{Rebecca R.G.}
\address{Department of Mathematical Sciences \\ George Mason University \\ Fairfax, VA  22030}
\email{rrebhuhn@gmu.edu}

\author{Hugh Geller}
\address{School of Mathematical and Data Sciences, West Virginia University, Morgantown, WV 26505}
\email{hugh.geller@mail.wvu.edu}

\title{Canonical forms of neural ideals}

\date{\today}
\begin{document}

\begin{abstract}
    Neural ideals, originally defined in \cite{neuralring}, give a way of translating information about the firing pattern of a set of neurons into a pseudomonomial ideal in a polynomial ring. We give a simple criterion for determining whether a neural ideal is in canonical form, along with an improved algorithm for computing the canonical form of a neural ideal.
\end{abstract}

\maketitle
\setcounter{tocdepth}{1} 

\section{Introduction}

Neural rings were first introduced in \cite{neuralring} as an algebraic tool to study \textit{receptive field codes} (RF codes). These are codes formed by measuring neuron activity in response to stimuli spaces. Typically, we consider stimuli spaces as subsets of a subspace $X$ of $\mathbb{R}^d$.   One of the motivating ideas in the field is to identify within a given code whether or not it arises from a stimuli space consisting of convex sets \cite{MR3595298,MR3633775,MR3914562,MR3944460,MR4042685,MR3954293,MR3998832}.

The neural ring serves as a bridge allowing for algebraic methods to be applied to neuroscience and coding theory.  Given $n$ neurons, the first step of the construction translates neuron firings to a binary code $\mathcal{C} \subseteq \{0,1\}^n$ where each possible combination of neuron firing is translated into a codeword of length $n$. For a given codeword $c = (c_1,\ldots,c_n)$, the entry $c_i = 1$ means the $i$th neuron is firing whereas $c_i = 0$ means the $i$th neuron is not firing. Given a neural code $\mathcal{C}$, the neural ring is $R/I_{\mathcal{C}}$, where $R=\F_2[x_1,\ldots,x_n]$ and $I_{\mathcal{C}}$ is the vanishing ideal of $\mathcal{C}$.
Using a vanishing ideal enables us to apply tools from algebraic geometry to neural codes.
However, from the perspective of distinguishing neural codes, this ideal is too large as it includes the trivial relations $\{x_i(1 - x_i): 1 \leq i \leq n\}$; these hold for every neural code and so are not useful in distinguishing neural codes from each other.
Instead, \cite{neuralring} defines the \textit{neural ideal} $J_{\mathcal{C}} \subset I_{\mathcal{C}} \subseteq R$ formed from the non-trivial relations of the neural code. The two ideals relate through the equation \[I_{\mathcal{C}} = J_{\mathcal{C}} + \langle x_i(1 - x_i) : 1 \leq i \leq n \rangle\] (cf. \cite[Lemma 3.2]{neuralring}). As such, one can study the neural ideal in place of the vanishing ideal.

The information encoded in the neural ideal is made more accessible by presenting the ideal in \textit{canonical form} (cf. \cite[Section 4.3]{neuralring}). The canonical form of a neural ideal can be used to identify various obstructions to convexity \cite{MR3944460}. As such it is useful to know whether or not a neural ideal is being presented in canonical form. 

One downside to the neural ideal from an algebraic perspective is that it is generated by \textit{pseudomonomials}, products $\Pi_{i \in \sigma} x_i \Pi_{i \in \tau} (1-x_i)$, which are neither graded nor local and as a result are not as well studied in commutative algebra. To deal with this issue, G\"{u}nt\"{u}rk\"{u}n, Jeffries, and Sun give a process for \textit{polarizing} a neural ideal, turning it into a squarefree monomial ideal in the extension ring $\F_2[x_1,\ldots,x_n,y_1,\ldots,y_n]$. 
In our work, we utilize the polarization technique from \cite{polarizationofneuralrings} to apply monomial ideal results to classify when certain (polarized) neural ideals are presented in canonical form. As a result, our main theorem below establishes sufficient conditions on the generators of a polarized neural ideal to determine when it is in canonical form. Here, saying that two monomials $g$ and $h$ ``share an index" means that there is some $1 \le i \le n$ such that $x_i \mid g$ and $y_i \mid h$ or vice versa.

\begin{thm}[Theorem \ref{thm:mainthm}]
Let $\mathfrak{a}=(g_1,\ldots,g_k)$ be a polarized neural ideal such that $g_{j_1} \nmid g_{j_2}$ for any $1 \le j_1 \ne j_2 \le k$, and $x_iy_i \nmid g_j$ for any $1 \le i \le n$ and $1 \le j \le k$. 
If for some pair $g_{j_1},g_{j_2}$ of generators of $\mathfrak{a}$, $g_{j_1}$ and $g_{j_2}$ share exactly 1 index $i$ and no other generator of $\mathfrak{a}$ divides $\frac{[g_{j_1}g_{j_2}]}{x_iy_i}$, then $\mathfrak{a}$ is not in canonical form.
Otherwise, $\mathfrak{a}$ is in canonical form.
\end{thm}

As a consequence of this theorem, we are able to give a shortened algorithm for computing the canonical form of a neural ideal (Algorithm \ref{alg:shorter}). We are also able to list several families of polarized neural ideals that are \underline{not} in canonical form and give the canonical form for each of these familes (see~Section \ref{sec:computations}). It is important to note that none of our results depend on the number of neurons being studied; that is, for $\mathcal{C} \subseteq \{0,1\}^n$, our results do not put any conditions on the value of $n$.

The structure of our paper is as follows. Section \ref{sec:background} covers the necessary background information from \cite{neuralring} and \cite{polarizationofneuralrings}. We also recall some useful observations about squarefree monomial ideals.

Section \ref{sec:algorithm} is dedicated to translating the algorithm from \cite{neuralring} for computing the canonical form of a neural ideal into an algorithm that can be applied to polarized neural ideals (Algorithm \ref{alg:canonform}). In Section \ref{sec:keylemmas}, we identify common patterns in the algorithm based off of the generators of our ideal, giving us a number of shortcuts through Algorithm \ref{alg:canonform}.

Along with containing our main theorem, Section \ref{sec:mainresults} includes a complete classification of the canonical forms of all two-generated neural ideals and a simplified algorithm for computing the canonical form of a neural ideal (Algorithm \ref{alg:shorter}).
In Section \ref{sec:genericcanonicalforms}, we discuss a \textit{generic} canonical form, a tool that allows us to use the canonical form of a single neural ideal to compute the canonical forms of a number of related neural ideals.
Section \ref{sec:computations} is dedicated to developing techniques for determining the canonical forms of families of neural ideals.

\section{Background on neural ideals and polarization}
\label{sec:background}

Throughout the paper, we make use of the following notation.

\begin{notation}
\phantom{making space}
\begin{enumerate}
    \item $R$ will denote the ring $\mathbb{F}_2[x_1,\ldots,x_n]$.
    \item A \textit{pseudomonomial} is a product $\prod_{i \in \sigma} x_i \prod_{i \in \tau} (1-x_i)$ in $R$ where $\sigma,\tau \subseteq \{1,2,\ldots,n\}$ and $\sigma \cap \tau = \emptyset$.
\end{enumerate}
\end{notation}

\begin{defn}
To get the neural ideal, we work with a space $X \subseteq \R^d$ and open sets $U_1,\ldots,U_n \subseteq X$ where neuron $i$ fires on the set $U_i$.
The \textit{neural ideal}, as introduced in \cite{neuralring}, is an ideal of $R$ that captures the \textit{RF-structure} of a neural code, specifically the relations
\[\bigcap_{i \in \sigma} U_i \subseteq \bigcup_{i \in \tau} U_i,\]
where $\sigma,\tau \subseteq \{1,\ldots,n\}$ and $\sigma \cap \tau= \emptyset$. More precisely, the neural ideal is generated by the pseudomonomials $\prod_{i \in \sigma} x_i \prod_{i \in \tau} (1-x_i)$ corresponding to the relations above \cite[Section 4.2]{neuralring}.
\end{defn}

Note that the Boolean relations $U_i \subseteq U_i$ always hold, but they do not distinguish between distinct neural codes, so they are not included in the neural ideal. Hence we assume that no generator of a neural ideal is divisible by $x_i(1-x_i)$ for any $1 \le i \le n$.

\begin{defn}
The \textit{canonical form} of the neural ideal captures the minimal relations in the RF-structure of the neural code, in the sense that the intersection on the left and the union on the right are irredundant (i.e., if we removed any $U_i$ from either side of the inclusion, the inclusion would no longer hold). It is generated by the pseudomonomials corresponding to the minimal relations \cite[Section 4.3]{neuralring}.
\end{defn}

One downside to these ideals is that they are generated by pseudomonomials, and as a result are neither graded nor local. In order to study objects like minimal free resolutions of these ideals, G\"{u}nt\"{u}rk\"{u}n, Jeffries, and Sun \cite{polarizationofneuralrings} developed a method to pass from pseudomonomial ideals in $R$ to squarefree monomial ideals in $S=\F_2[x_1,\ldots,x_n,y_1,\ldots,y_n]$ by \textit{polarization}.

\begin{defn}[\cite{polarizationofneuralrings}]
\phantom{creating space}
\begin{enumerate}
    \item $S$ will denote the extension ring $\mathbb{F}_2[x_1,\ldots,x_n,y_1,\ldots,y_n]$ of $R$.
    \item For a monomial $m$ we denote its largest square free divisor as $[m]$. For example, $[x^3y^2] = xy$. One can quickly check that if $m_1$ and $m_2$ are squarefree monomials, then $[m_1m_2] = \text{lcm}\{m_1,m_2\}$.
    \item We will refer to a \textit{depolarization} map $d:S \to R$, which is the ring homomorphism sending $x_i \mapsto x_i, y_i \mapsto 1-x_i$. The map $d$ induces an isomorphism of $S$-modules $R \cong S/\left(\sum_{i=1}^n (x_i+y_i-1)\right)$. \cite{polarizationofneuralrings}
    \item We work with neural ideals, which are generated by pseudomonomials.
    \item We will also refer to a \textit{polarization} function $\cP$ (not a homomorphism) sending pseudomonomials in $R$ to squarefree monomials in $S$ via \[\prod_{i \in \sigma} x_i \prod_{i \in \tau} (1-x_i) \mapsto \prod_{i \in \sigma} x_i \prod_{i \in \tau} y_i.\]
    \item If $\mathfrak{a}$ is an ideal of $R$, then $\cP(\mathfrak{a})$ is the ideal of $S$ generated by $\cP(f)$ for all $f \in \mathfrak{a}$.
\end{enumerate}
\end{defn}

We will often work with the polarization of a neural ideal, while explaining how this translates to the depolarized version.

\begin{thm}[{\cite[Theorem 3.2]{polarizationofneuralrings}}]
\label{thm:canonpolarizes}
Let $\cA=(g_1,\ldots,g_k) \subseteq R$ be a neural ideal in canonical form. Then $\cP(\cA)=(\cP(g_1),\ldots,\cP(g_k))$.
\end{thm}

 However, this does not always hold when the ideal is not in canonical~form. 
 
 \begin{example}
 For example, $(x_1,x_2(1-x_1))$ contains $x_2$, so $\cP(x_1,x_2(1-x_1))$ contains $x_2$. However, \[x_2 \not\in (\cP(x_1),\cP(x_2(1-x_1)))=(x_1,x_2y_1).\] The canonical form of $(x_1,x_2(1-x_1))$ turns out to be $(x_1,x_2)$.
 \end{example}

\begin{defn}[\cite{neuralring}]
A \textit{pseudomonomial prime} of $R$ is an ideal $p$ generated by a subset of $\{x_1,\ldots,x_n,1-x_1,\ldots,1-x_n\}$ such that for each $1 \le i \le n$, $p$ does not contain both $x_i$ and $1-x_i$.
\end{defn}

\begin{thm}[{\cite[Theorem 5.4]{neuralring}}]
\label{thm:pseudomonomialprimarydecompunique}
Let $\mathfrak{a}$ be a pseudomonomial ideal in $\mathbb{F}_2[x_1,\ldots,x_n]$. Then the primary decomposition of $\mathfrak{a}$ consists of pseudomonomial primes, and so $\mathfrak{a}$ has no embedded primes. Hence its primary decomposition is unique.
\end{thm}

\begin{thm}[{\cite[Part 3]{monomialideals}}]
\label{thm:squarefreemonomialprimary}
Let $\mathfrak{a}$ be a squarefree monomial ideal in a polynomial ring $k[x_1,\ldots,x_n]$. The primary decomposition of $\mathfrak{a}$ consists of primes generated by a subset of $\{x_1,\ldots,x_n\}$ (\textit{monomial primes}), and so $\mathfrak{a}$ has no embedded primes. Hence its primary decomposition is unique.
\end{thm}

\begin{rem}
Let $g_1,\ldots,g_k$ be squarefree monomials in a polynomial ring $k[x_1,\ldots,x_n]$. If $g_j=hh'$ for any squarefree monomials $h,h'$ in the ring, then
\[(g_1,\ldots,g_j,\ldots,g_k)=(g_1,\ldots,h,\ldots,g_k) \cap (g_1,\ldots,h',\ldots,g_k).\]
In order to compute the primary decomposition of $(g_1,\ldots,g_k)$, we may use this rule repeatedly until we get an intersection of monomial primes.
\end{rem}

\section{Computing the canonical form of a polarized neural ideal}
\label{sec:algorithm}

In this section we give a polarized version of the algorithm for computing the canonical form of a neural ideal and prove that it agrees with the original algorithm in \cite[Section 4.5]{neuralring}.

We recall the original algorithm for computing the canonical form of a neural ideal below.

\begin{algorithm}[{\cite[Section 4.5]{neuralring}}]
\label{alg:orig}
\phantom{making space}
\begin{enumerate}
    \item Start with a neural ideal $\mathcal{A}=(g_1,\ldots,g_k)$ in $R$.
    \item Compute the primary decomposition of $\mathcal{A}$. By Theorem \ref{thm:pseudomonomialprimarydecompunique}, the ideals $p_1,\ldots,p_s$ in the primary decomposition will all be generated by a subset of $\{x_1,\ldots,x_n,1-x_1,\ldots,1-x_n\}$.
    \item Compute the set of products $h_1 \cdots h_s$ where $h_\ell$ is a generator of $p_\ell$.
    \item Set $x_i(1-x_i)=0$ for each $1 \le i \le n$, and as a result remove any product divisible by $x_i(1-x_i)$ and replace any power $x_i^t$ by $x_i$ and $(1-x_i)^t$ by $1-x_i$. All remaining elements are now pseudomonomials.
    \item Remove any product that is a multiple of a product of lower degree.
\end{enumerate}
The remaining products give the canonical form of $\mathcal{A}$ in $R$.
\end{algorithm}

We give our polarized version of the algorithm below. 

\begin{algorithm}
\label{alg:canonform}
\phantom{making space}
\begin{enumerate}
    \item Start with a squarefree monomial ideal $\mathfrak{a}=(g_1,\ldots,g_k) \in S$. 
    \item Compute the primary decomposition of $\mathfrak{a}$. By Theorem \ref{thm:squarefreemonomialprimary}, the ideals $p_1,\ldots,p_s$ in the primary decomposition will all be generated by a subset of $\{x_1,\ldots,x_n,y_1,\ldots,y_n\}$.
    \item Set $x_i+y_i=1$ for $1 \le i \le n$, and as a result remove any ideal in the primary decomposition containing both $x_i$ and $y_i$ for any $1 \le i \le n$.
    \item Compute the intersection of the remaining ideals. Since the $p_{\ell}$ are squarefree monomial ideals, this is equivalent to computing $[h_1 \cdots h_s]$ for every set of choices of $h_{\ell}$ a generator of $p_{\ell}$.
    \item Impose the relations $x_iy_i=0$ for $1 \le i \le n$, and as a result remove any generator divisible by $x_iy_i$ for some $i$. 
    (These come from imposing the relations $x_i(1-x_i)=0$ in the depolarized ideal.) 
    \item Remove any generator that is a multiple of another generator.
\end{enumerate}
The remaining generators give the canonical form of $\mathfrak{a}$ in $S$.
\end{algorithm}

\begin{lemma}
\label{lem:polarization}
\phantom{making space}
\begin{enumerate}
    \item Let $f=\prod_{i \in \sigma} x_i \prod_{i \in \tau} (1-x_i)$ be a pseudomonomial in $R$. Then $d(\cP(f))=f$.
    \item Let $g=\prod_{i \in \sigma} x_i \prod_{i \in \tau} y_i$ be a monomial in $S$. Then $\cP(d(g))=g$.
    \item Let $\mathfrak{a}=(g_1,\ldots,g_{k})$ be a squarefree monomial ideal in $S$. Then $d(\mathfrak{a})~=~(d(g_1),\ldots,d(g_{k}))R$, though the generators may not be minimal. If $x_i,y_i \in \mathfrak{a}$ for some $1 \le i \le n$, then $d(\mathfrak{a})=(1)$.
    \item Let $p$ be a pseudomonomial prime of $R$. Then $p$ is in canonical form.
    \item Let $p$ be a pseudomonomial prime of $R$. Then $\cP(p)$ is a monomial prime of $S$.
    \item Let $\mathfrak{a}$ be a squarefree monomial ideal of $S$. If for all $1 \leq i \leq n$ $x_iy_i$ does not divide any generator, then $d(\mathfrak{a})$ is an ideal of $R$ generated by pseudomonomials. In particular, if $p$ is a monomial prime of $S$ such that for all $1 \le i \le n$, $(x_i,y_i)$ is not a subideal of $p$, then $d(p)$ is a pseudomonomial prime of $R$.
    \item Let $\cA$ be a pseudomonomial ideal of $R$. Then $\cA \subseteq d(\cP(\cA))$. If $\cA$ is a pseudomonomial prime, this is an equality.
    \item Let $p$ be a monomial prime of $S$ such that for $1 \le i \le n$, at most one of $x_i,y_i$ is in $p$. Then $\cP(d(p))=p$.
\end{enumerate}
\end{lemma}

\begin{proof}
Parts (1) and (2) are clear. Part (3) follows from observing that $d$ is an $S$-module homomorphism.

For (4), we note that $p$ is in canonical form since it is preserved by Algorithm \ref{alg:orig}. Then (5) follows by Theorem \ref{thm:canonpolarizes}.

Statement (6) follows from (3) and the fact that an ideal generated by some $x_i$ and some $1-x_j$ where there is no index for which we have both $x_i$ and $1-x_i$ as generators is prime.

To prove (7), first let $f \in \cA$ be a pseudomonomial. Then by part (1), $f=d(\cP(f)) \in d(\cP(\cA))$. If $\cA$ is a pseudomonomial prime, then without loss of generality it has the form $(x_1,\ldots,x_k,(1-x_{k+1}),\ldots,(1-x_{k+t}))$. Since this ideal is in canonical form by (4), it follows from Theorem \ref{thm:canonpolarizes} that $\cP(\cA)=(x_1,\ldots,x_k,y_{k+1},\ldots,y_{k+t})$. Then $d(\cP(\cA))=\cA$.

For (8), assume without loss of generality that \[p=(x_1,\ldots,x_k,y_{k+1},\ldots,y_{k+t}).\] Then $d(p)=(x_1,\ldots,x_k,(1-x_{k+1}),\ldots,(1-x_{k+t}))$ is in canonical form by (4). The result now follows from Theorem \ref{thm:canonpolarizes}.
\end{proof}

\begin{example}
Consider the ideal $\mathfrak{a}=(x_1,x_2y_1)$ of $S$. The depolarized ideal is $d(\mathfrak{a})=(x_1,x_2(1-x_1))=(x_1,x_2)$, so we see that there is a ``more minimal" way to write $d(\mathfrak{a})$, as mentioned in part (3) of Lemma \ref{lem:polarization}.
\end{example}

\begin{lemma}
\label{lem:primarydecomp}
Let $\cA \subseteq R$ be a neural ideal, and let $\mathfrak{a} \subseteq S$ be its polarization. The set of associated primes of $\cA$ is equal to the set of primes obtained by depolarizing the set of associated primes of $\mathfrak{a}$, after removing primes containing both $x_i$ and $y_i$ for some $1 \le i \le n$.
\end{lemma}

\begin{proof}
We need to prove two claims:
\begin{enumerate}
    \item If $p$ is a minimal prime of $\cA$, then $\cP(p)$ is a minimal prime of $\mathfrak{a}$.
    \item If $p$ is a minimal prime of $\mathfrak{a}$ and there is no $1 \le i \le n$ such that $p$ contains both $x_i$ and $y_i$, then $d(p)$ is a minimal prime of $\cA$.
\end{enumerate}
Note that by Theorem \ref{thm:squarefreemonomialprimary},
the irredundant primary decomposition of $\mathfrak{a}$ is unique. Similarly, by Theorem \ref{thm:pseudomonomialprimarydecompunique}, the irredundant primary decomposition of $\cA$ is also unique.

\begin{proof}[Proof of Claim 1]
Let $p$ be a minimal prime of $\cA$. By Theorem \ref{thm:pseudomonomialprimarydecompunique}, $p$ is generated by pseudomonomials $x_i$ and $1-x_i$ (where for each $i$, at most one of $x_i,1-x_i$ is a generator). By Lemma \ref{lem:polarization}, $\cP(p)$ is prime, and $\mathfrak{a} \subseteq \cP(p)$. Suppose $\cP(p)$ is not minimal, i.e. there is a prime $q$ with $\mathfrak{a} \subseteq q \subset \cP(p)$. Since $\mathfrak{a}$ is a squarefree monomial ideal, $q$ is generated by a subset of the variables $x_1,\ldots,x_n,y_1,\ldots,y_n$.
Since $q \subseteq \cP(p)$, which does not contain both $x_i$ and $y_i$ for any $i$, Lemma \ref{lem:polarization} implies $d(q)$ is prime in $R$. We have
\[\cA \subseteq d(\cP(\cA))=d(\mathfrak{a}) \subseteq d(q) \subseteq d(\cP(p))=p.\]
Since $p$ was a minimal prime of $\cA$ by assumption, $d(q)=p$. Suppose that $q \subsetneqq \cP(p)$. Then there must be some $z \in \{x_1,\ldots,x_n,y_1,\ldots,y_n\}$ such that $z \in \cP(p)$, but $z \not\in q$. Then $d(z) \in d(\cP(p))=p=d(q)$. Applying Lemma \ref{lem:polarization}, $z=\cP(d(z)) \in \cP(d(q))=q$, since $q$ does not contain $x_i$ and $y_i$ for any $1 \le i \le n$.
So $q=\cP(p)$, which implies that $\cP(p)$ is a minimal prime of $\mathfrak{a}$.
\end{proof}

\begin{proof}[Proof of Claim 2]
Let $p$ be a minimal prime of $\mathfrak{a}$ such that for $1 \le i \le n$, $p$ does not contain both $x_i$ and $y_i$. Then by Lemma \ref{lem:polarization}, $d(p)$ is a pseudomonomial prime and 
\[\cA \subseteq d(\cP(\cA))=d(\mathfrak{a}) \subseteq d(p).\]
Suppose there is some pseudomonomial prime $q$ of $R$ such that $\cA \subseteq q \subset d(p)$. Then $\cP(q)$ is prime and $\mathfrak{a} \subseteq \cP(q) \subseteq \cP(d(p))$, which is equal to $p$ by Lemma \ref{lem:polarization}. Since $p$ is a minimal prime of $\mathfrak{a}$, $\cP(q)=p$. Then by Lemma \ref{lem:polarization} and since $d$ is a homomorphism, $q=d(\cP(q))=d(p)$. So $d(p)$ is a minimal prime of $\cA$, as desired.
\end{proof} 
Since neither $\cA$ nor $\mathfrak{a}$ has any embedded primes, this proves the result.
\end{proof}

\begin{prop}
\label{pr:algworks}
Let $\cA \subseteq \mathbb{F}_2[x_1,\ldots,x_n]$ be a neural ideal, and let $\mathfrak{a}$ be its polarization. The canonical form of $\cA$ as computed in Algorithm \ref{alg:orig} is equal to the depolarization of the canonical form of $\mathfrak{a}$ as computed in Algorithm~\ref{alg:canonform}. 
\end{prop}

\begin{proof}
First, by Lemma \ref{lem:primarydecomp}, at the end of Step 3 of Algorithm \ref{alg:canonform}, the remaining set of prime ideals, when depolarized, is equal to the set of prime ideals occurring in the primary decomposition of $\cA$. 

Next we prove that the set of products left at the end of Step 5 of Algorithm \ref{alg:canonform}, depolarized, agree with the products left at the end of Step 4 of Algorithm \ref{alg:orig}. We need to show that if $z$ and $w$ are pseudomonomials in $R$, if we apply the relations $x_i(1-x_i)=0$ to $zw$, the result is the same as if we compute $[\cP(z)\cP(w)]$, impose the relation $x_iy_i=0$, and then depolarize. Suppose
\[
z=\prod_{i \in \sigma_1} x_i \prod_{i \in \tau_1} (1-x_i), \quad \text{and} \quad  w=\prod_{i \in \sigma_2} x_i \prod_{i \in \tau_2} (1-x_i).
\]
Then the result of applying the first process to $zw$ is
\[
\begin{cases}
\prod_{i \in (\sigma_1 \cup \sigma_2)} x_i \prod_{i \in (\tau_1 \cup \tau_2)} (1-x_i) & (\sigma_1 \cup \sigma_2) \cap (\tau_1 \cup \tau_2)=\emptyset \\
0 & \text{else}.\\
\end{cases}
\]
The result of applying the second process to $[\cP(z)\cP(w)]$ is
\[
\begin{cases}
\prod_{i \in (\sigma_1 \cup \sigma_2)} x_i \prod_{i \in (\tau_1 \cup \tau_2)} y_i & (\sigma_1 \cup \sigma_2) \cap (\tau_1 \cup \tau_2)=\emptyset \\
0 & \text{else}.\\
\end{cases}
\]
Depolarizing, this agrees with the result of the first process.

Finally, by \cite[Lemma 3.1]{polarizationofneuralrings}, if $z$ and $w$ are two pseudomonomials in $R$, $z \mid w$ if and only if $\cP(z) \mid \cP(w)$. Hence the depolarization of the result of Step 6 of Algorithm \ref{alg:canonform} agrees with the result of Step 5 of Algorithm~\ref{alg:orig}.
\end{proof}

\section{Deconstructing the canonical form algorithm}
\label{sec:keylemmas}

In this section we prove a number of lemmas describing in more detail what happens to the generators of a polarized neural ideal under the steps of Algorithm \ref{alg:canonform} and giving us shortcuts through the algorithm.

\begin{defn}
\label{def:sharedindex}
We say two monomials $m_1$ and $m_2$ \textit{share an index} $i$ if $x_i$ divides one of them and $y_i$ divides the other.
If a single monomial is divisible by $x_iy_i$, we do not count this as sharing an index with itself (see Remark \ref{rem:sharedindeximportance}).
Moreover, we say that the generators of a squarefree monomial ideal $\mathfrak{a}$ \textit{share the index} $i$ if there are two generators of the ideal sharing the index $i$.
\end{defn}

\begin{example}
For example, $x_1y_2$ and $x_3y_1$ share the index 1. However, $x_1y_2$ and $x_1y_3$ share no index, even though $x_1$ divides both.
\end{example}

\begin{notation}
We refer to the generators of an ideal as $g_1,\ldots,g_k$, frequently pulling out $x_i$ and $y_i$ when $i$ is a shared index between the $g_j$. For example, we write $(x_1g_1,y_1g_2)$. The $g_j$ will always be squarefree monomials.

For ease of notation, we write $x_1,x_2,\ldots$ and $y_1,y_2,\ldots$ instead of $x_{i_1},x_{i_2},\ldots$ and $y_{i_1},y_{i_2},\ldots$. However, all of our results hold if we permute the indices $1,\ldots,n$ on the variables.
\end{notation}

\begin{rem}
\label{rem:squarefreereason}
We state the results in this section for squarefree monomial ideals rather than for polarized neural ideals so that we can apply them to ideals appearing as intermediate stages in Algorithm \ref{alg:canonform}, which may have generators divisible by $x_iy_i$ for some $1 \le i \le n$.
\end{rem}

\begin{defn}
Let $\mathfrak{a}$ be a squarefree monomial ideal in $S$. We say that $\mathfrak{a}$ is \textit{\fourstep} if it is the result of applying Steps 1-4 of Algorithm \ref{alg:canonform} to some squarefree monomial ideal, and call process of applying Steps 1-4 to a squarefree monomial ideal \textit{\fourstepping} it. We say that $\mathfrak{a}$ is in \textit{\fivestep} if it is the result of applying Steps 1-5 of Algorithm \ref{alg:canonform} to some squarefree monomial ideal.
\end{defn}

\begin{lemma}
\label{lem:oneatatime}
Let $\mathfrak{a}$ be a squarefree monomial ideal in $S$. For some $0 \le t \le n$ and $i_1,\ldots,i_t \in \{1,\ldots,n\}$, we can write
$\mathfrak{a}=\mathfrak{a}_0 \cap \mathfrak{a}_1 \cap \ldots \cap \mathfrak{a}_t,$
where 
\begin{enumerate}
    \item $\mathfrak{a}_0$ is an intersection of monomial primes, none of which contain both $x_i$ and $y_i$ for any $1 \le i \le n$, and
    \item for $1 \le j \le t$, $\mathfrak{a}_j$ is an intersection of monomial prime ideals containing both $x_{i_j}$ and $y_{i_j}$.
\end{enumerate}
For any such decomposition, \fourstepping\ $\mathfrak{a}$ returns $\mathfrak{a}_0$, and so $\mathfrak{a}$ and $\mathfrak{a}_0$ have the same canonical form. In fact for any $s \le t$, \fourstepping\
$\mathfrak{a}'=\mathfrak{a}_0 \cap \mathfrak{a}_s \cap \ldots \cap \mathfrak{a}_t$
yields $\mathfrak{a}_0$,
and so $\mathfrak{a}'$ has the same canonical form as $\mathfrak{a}_0$.

Consequently, if $\mathfrak{a}=\mathfrak{b} \cap \mathfrak{c}$ and $(x_i,y_i) \subseteq \mathfrak{b}$ for some $1 \le i \le n$, recomposing $\mathfrak{a}$ returns the same result as recomposing $\mathfrak{c}$.
\end{lemma}

\begin{proof}
We prove this by induction on $t$, the order of the set
\[I=\{i:1 \le i \le n,(x_i,y_i) \subseteq p \text{ for some } p \text{ in a primary decomposition of } \mathfrak{a}\}.\]
If $t=0$, $\mathfrak{a}=\mathfrak{a}_0$ so the result is immediate. Suppose $t>0$. Without loss of generality, assume that $I=\{1,\ldots,t\}$. For $1 \le i \le t$, set $\mathfrak{a}_i$ to be the intersection of the prime ideals $p$ appearing in a primary decomposition of $\mathfrak{a}$ that contain $x_i$ and $y_i$.
Set $\mathfrak{a}_0$ to be the intersection of the primes $p$ appearing in the primary decomposition of $\mathfrak{a}$ that do not contain $x_i,y_i$ for any $1 \le i \le n$.

When we perform Step 3 of Algorithm \ref{alg:canonform}, we will remove exactly the primes $p$ that contain a pair $x_i,y_i$ for some $1 \le i \le t$. This exactly corresponds to removing $\mathfrak{a}_1,\ldots,\mathfrak{a}_t$, leaving us with $\mathfrak{a}_0$ when we recompose $\mathfrak{a}$ in Step 4.

The same argument holds for $\mathfrak{a}'$.

For the last statement, if $\mathfrak{a}=\mathfrak{b} \cap \mathfrak{c}$, where $(x_i,y_i) \subseteq \mathfrak{b}$ for some $1 \le i \le n$, every minimal prime of $\mathfrak{b}$ will be part of $\mathfrak{a}_s$ for some $1 \le s \le t$, so they will all be removed in recomposing $\mathfrak{a}$. Thus it suffices to recompose $\mathfrak{c}$.
\end{proof}

\begin{defn}
 We call applying Steps 1-4 of Algorithm \ref{alg:canonform} to a squarefree monomial ideal $\mathfrak{a} \subseteq S$ but only removing primes containing $(x_i,y_i)$ for a particular $i$ \textit{recomposing $\mathfrak{a}$ with respect to the index $i$}. We may also refer to recomposing $\mathfrak{a}$ with respect to a subset of $\{1,\ldots,n\}$, by which we mean we recompose with respect to each index in the set. If no set is specified, we are recomposing with respect to $\{1,\ldots,n\}$.
\end{defn}

\begin{rem}
\label{rem:sharedindeximportance}
We observe that a minimal prime of $\mathfrak{a}$ contains $(x_i,y_i)$ if and only if at least two generators of $\mathfrak{a}$ share the index $i$. The backward direction follows from the process of decomposing a squarefree monomial ideal. For example, $(x_1,x_2y_1)=(x_1,x_2) \cap (x_1,y_1)$.
For the forward direction, if some minimal prime of $\mathfrak{a}$ contains $(x_i,y_i)$, then $\mathfrak{a}$ must contain at least one multiple of $x_i$ and at least one (distinct) multiple of $y_i$.

As a result, the indices $i_1,\ldots,i_t$ from Lemma \ref{lem:oneatatime} are the indices shared by the generators of $\mathfrak{a}$. Hence it suffices to recompose an ideal with respect to the indices shared by its generators, as recomposing with respect to a non-shared index removes no minimal primes and thus preserves the ideal.

This also explains why, in Definition \ref{def:sharedindex}, we do not count $i$ as a shared index among the generators of an ideal $\mathfrak{a}$ if $x_iy_i$ divides a single generator, but neither $x_i$ nor $y_i$ divides any other generator of $\mathfrak{a}$ (see Definition \ref{def:sharedindex}): such a generator does not give us a minimal prime containing $(x_i,y_i)$. For example, $(x_1y_1,x_2)=(x_1,x_2) \cap (y_1,x_2)$.
\end{rem}

\begin{rem}
\label{rem:removemultiples}
In computing canonical forms, we may remove generators that are a multiple of another generator (or duplicate generators) at any stage of the algorithm, rather than waiting until Step 6. Removing such generators preserves the ideal, and hence the primary decomposition as well as the intersection that occurs in Step 4 of Algorithm \ref{alg:canonform}.
\end{rem}

\begin{lemma}
\label{lem:combinecanonical}
Let $t \ge 2$. Let 
\[\mathfrak{a}=\bigcap_{j=1}^t \mathfrak{a}_j,\]
where $\mathfrak{a},\mathfrak{a}_j$ are squarefree monomial ideals of $S$. Recomposing $\mathfrak{a}$ with respect to an index $i$ is equivalent to recomposing the $\mathfrak{a}_j$ with respect to the index $i$ and then intersecting the results.
\end{lemma}

\begin{proof}
We prove the case $t=2$. The rest follows by induction. We show that the set of primes $\{p_i\}$ appearing in a primary decomposition of $\mathfrak{a}$ is equal to the subset of the primes $\{q_i^{(1)},q_j^{(2)}\}$ appearing in the primary decompositions of $\mathfrak{a}_1$ and $\mathfrak{a}_2$ obtained by removing duplicates and primes that are now non-minimal: 
 we first prove that any prime minimal over $\mathfrak{a}$ is minimal over $\mathfrak{a}_1$ or $\mathfrak{a}_2$. Suppose $\mathfrak{p}$ is minimal over $\mathfrak{a}$ but $\mathfrak{a}_1 \not\subseteq \mathfrak{p}$; we will show $\mathfrak{p}$ is minimal over $\mathfrak{a}_2$. Since $\mathfrak{a}_1 \not\subseteq \mathfrak{p}$, there exists $\alpha \in \mathfrak{a}_1 \backslash \mathfrak{p}$. We observe that $\alpha \mathfrak{a}_2 \subseteq \mathfrak{a}_1 \cap \mathfrak{a}_2 = \mathfrak{a} \subseteq \mathfrak{p}.$ Since $\mathfrak{p}$ is prime and $\alpha \notin \mathfrak{p}$, we must have  $\mathfrak{a}_2 \subseteq \mathfrak{p}$.

To see that $\mathfrak{p}$ is minimal over $\mathfrak{a}_2$, we consider another prime $\mathfrak{p}'$ such that $\mathfrak{a}_2 \subseteq \mathfrak{p}' \subseteq \mathfrak{p}$. Since $\mathfrak{a} \subseteq \mathfrak{a}_2$, we have $\mathfrak{a} \subseteq \mathfrak{p}' \subseteq \mathfrak{p}$. Since $\mathfrak{p}'$ is prime and $\mathfrak{p}$ is a minimal prime over $\mathfrak{a}$, we must have $\mathfrak{p}' = \mathfrak{p}$.

Next, we consider a prime $\mathfrak{p}$ that is minimal over $\mathfrak{a}_1$ but not over $\mathfrak{a}$. In that case, there is another prime $\mathfrak{p}'$ minimal over $\mathfrak{a}$ such that $\mathfrak{a} \subseteq \mathfrak{p}' \subsetneq \mathfrak{p}$. By the minimality of $\mathfrak{p}$ over $\mathfrak{a}_1$, we know that $\mathfrak{a}_1 \not\subseteq \mathfrak{p}'$. Since $\mathfrak{p}'$ is minimal over $\mathfrak{a}$ but not over $\mathfrak{a}_1$, we find that $\mathfrak{p}'$ is minimal over $\mathfrak{a}_2$. Consequently, if we take all the primes minimal over $\mathfrak{a}_1$ or $\mathfrak{a}_2$ and remove any primes that are not minimal over $\mathfrak{a}$ or are duplicates, we obtain all the primes minimal over $\mathfrak{a}$.

This remains true after removing any primes containing $(x_i,y_i)$ for some $1 \le i \le n$ from both sets. Computing the intersection of the remaining $p_i$ is equivalent to separately intersecting the remaining $q_i^{(1)}$ and $q_i^{(2)}$, and then intersecting the two resulting ideals.
\end{proof}

\begin{example}
We use Lemmas \ref{lem:oneatatime} and \ref{lem:combinecanonical} to compute the canonical form of $\mathfrak{a}=(x_1x_2x_3,x_4y_1y_2)$. We begin decomposing $\mathfrak{a}$ as follows:
\begin{align*}
    (x_1x_2x_3,x_4y_1y_2) &=(x_1,x_4y_1y_2) \cap (x_2,x_4y_1y_2) \cap (x_3,x_4y_1y_2)\\
    &=(x_1,y_1) \cap (x_1,x_4y_2) \cap (x_2,y_2) \cap (x_2,x_4y_1) \cap (x_3,y_4y_1y_2).
\end{align*}
We remove the 1st and 3rd pieces. Notice that it is impossible for any of the other 3 pieces to decompose into primes containing a pair $x_i,y_i$. By Lemma \ref{lem:oneatatime}, these pieces are recomposed, and by Lemma \ref{lem:combinecanonical}, their intersection is also recomposed,
so we can jump straight to Steps 4-6, performed together:
\begin{align*}
    (x_1,x_4y_2) \cap (x_2,x_4y_1) \cap (x_3,x_4y_1y_2) &\to (x_1x_2x_3,x_4y_1y_2).
\end{align*}
\end{example}

\begin{lemma}
\label{lem:decomp}
Let $\mathfrak{a}=(g_1,\ldots,g_k,h_1,\ldots,h_\ell)$ be a squarefree monomial ideal in $S$ such that for every $1 \le j_1 \le k$ and $1 \le j_2 \le \ell$, $g_{j_1}$ and $h_{j_2}$ do not share any indices (resp. the index $i$). Then recomposing $\mathfrak{a}$ (resp. with respect to the index $i$) gives the same result as recomposing  $(g_1,\ldots,g_k)$ and $(h_1,\ldots,h_k)$ (resp. with respect to the index $i$) and taking their sum.

More generally, recomposing the neural ideal \[\mathfrak{a}=(g_{1,1},\ldots,g_{1,k_1},g_{2,1},\ldots,g_{2,k_2},\ldots,g_{m,k_m})\] (resp. with respect to the index $i$),
where $g_{j_1,t_1}$ and $g_{j_2,t_2}$ do not share any indices (resp. the index $i$) for any $j_1 \ne j_2$,
agrees with recomposing each of the $(g_{j,1},\ldots,g_{j,k_j})$ (resp. with respect to the index $i$) and summing the results. 

Consequently, the canonical form (resp. \fivestep) of $\mathfrak{a}$ is the sum of the canonical forms (resp. \fivestep) of the $(g_{j,1},\ldots,g_{j,t_j})$, up to removing any generator divisible by another.
\end{lemma}

\begin{proof}
We prove the result in the first paragraph; the second paragraph follows by induction. We first prove the case where the $g_{j_1}$ and $h_{j_2}$ do not share the index $i$. Suppose the generators of $\mathfrak{a}$ share the index $i$. Without loss of generality, assume $i$ is a shared index for $(g_1,\ldots,g_k)$. Then by our hypotheses, $x_i$ and $y_i$ do not divide any $h_j$. Hence no minimal prime of $(h_1,\ldots,h_\ell)$ contains either $x_i$ or $y_i$.

Consequently, we can write 
$\mathfrak{a}=(G' \cap G_i)+H=(G'+H)\cap (G_i+H),$
where $H=(h_1,\ldots,h_{\ell})$ and $G',G_i$ are the components of $G=(g_1,\ldots,g_k)$ from the statement of Lemma \ref{lem:oneatatime}, where $G_i$ is the intersection of all minimal primes of $G$ containing $(x_i,y_i)$, and $G'$ is the intersection of the minimal primes of $G$ not containing $(x_i,y_i)$. In particular, $G'$ is equal to the result of \fourstepping\ $G$ with respect to the index $i$.
 The minimal primes of $G' + H$ are all of the form $\mathfrak{p}' + \mathfrak{p}$ where $\mathfrak{p}'$ is minimal over $G'$ and $\mathfrak{p}$ is minimal over $H$.
Since no minimal prime of $H$ contains $(x_i)$ or $(y_i)$, $G'+H$ is equal to the intersection of the minimal primes of $\mathfrak{a}$ not containing $(x_i,y_i)$, and $G_i+H$ is equal to the intersection of the minimal primes of $\mathfrak{a}$ containing $(x_i,y_i)$.

By Lemma \ref{lem:oneatatime} applied to $\mathfrak{a}$, \fourstepping\ $\mathfrak{a}$ is equal to \fourstepping\ $G'+H$ (with respect to any index or set of indices). Further, $G'+H$ has no minimal primes containing $(x_i,y_i)$ so it is already \fourstep\ with respect to the index $i$.
Finally, by Lemma \ref{lem:oneatatime}, \fourstepping\ $G'$ agrees with \fourstepping\ $(g_1,\ldots,g_k)$. As a result, if instead we \fourstep\ $G$ and $H$ separately with respect to the index $i$ and add the results, we will get $G'+H$. This proves the result for a single index. 

The multiple index result follows by recomposing $\mathfrak{a}$ one index at a time as in Lemma \ref{lem:oneatatime}; at each index $i$, no $g_{j_1}$ and $h_{j_2}$ share the index $i$.
The final paragraph of the result follows by computing the \fivestep\ or canonical form as appropriate.
\end{proof}

\begin{example}
The canonical form of $\mathfrak{a}=(x_1x_2,x_2y_1,x_3x_4y_5,x_2x_5y_4)$ in $S$ is equal to the sum of the canonical forms of the ideals $(x_1x_2,x_2y_1)$ and $(x_3x_4y_5,x_2x_5y_4)$, up to removing generators that are a multiple of another generator. Notice that both pieces contain $x_2$, but no monomial is divisible by $y_2$ so 2 is not a shared index.
\end{example}

\begin{lemma}
\label{lem:boringgens}
Let $\mathfrak{a}=(g_1,\ldots,g_k)$ be a squarefree monomial ideal of $S$. If for some $0 \le \ell \le k$ and every $\ell+1 \le j \le k$, $g_j$ shares no indices (resp. does not share the index $i$) with any other generator, then $g_{\ell+1},\ldots,g_k$ are generators of the recomposition of $\mathfrak{a}$ with respect to every $1 \le i \le n$ (resp. the index $i$). If in addition $x_iy_i \nmid g_j$ for any $\ell+1 \le j \le k$, they are generators of the \fivestep\ of $\mathfrak{a}$. The canonical form of $\mathfrak{a}$ can then be found by removing generators that are a multiple of another generator.

In particular, if $\mathfrak{a}$ does not have any shared indices, then the ideal is \fourstep. If in addition no generator is divisible by any $x_iy_i$ and no generator is a multiple of another, the ideal is in canonical form.
\end{lemma}

\begin{proof}
This follows from Lemma \ref{lem:decomp}.
\end{proof}

\begin{example}
By Lemma~\ref{lem:boringgens}, the neural ideal $(x_1y_2,x_3y_2,x_1x_4)$ is in canonical form.
\end{example}

\begin{lemma}
\label{lem:lotsofsamexiyi}
Let $\mathfrak{a}=(x_1g_{11},\ldots,x_1g_{1k_1},y_1g_{21},\ldots,y_1g_{2k_2},g_{31},\ldots,g_{3k_3}),$
be a squarefree monomial ideal in $S$,
where $k_1,k_2 \ge 0$ and the only factors of $x_1,y_1$ appearing among the generators of $\mathfrak{a}$ are the ones shown. Then recomposing $\mathfrak{a}$ with respect to the index 1 returns
\begin{align*}
    \mathfrak{a}' =& (x_1g_{11},\ldots,x_1g_{1k_1}, y_1g_{21},\ldots,y_1g_{2k_2}, [g_{1j_1}g_{2j_2}]_{1 \le j_1 \le k_1, 1 \le j_2 \le k_2}, g_{31},\ldots,g_{3k_3}).
\end{align*}
If the only indices shared by the generators of $\mathfrak{a}$ are the $x_1,y_1$ shown, $\mathfrak{a}'$ is the \fivestep\ of $\mathfrak{a}$. In this case the canonical form of $\mathfrak{a}$ is found by removing generators of $\mathfrak{a}'$ that are a multiple of another generator.

For any pair $g_{1j_1},g_{2j_2}$ that share an index, $[g_{1j_1}g_{2j_2}]$ will be removed in Step 5 of Algorithm \ref{alg:canonform}.
\end{lemma}

\begin{proof}
By Lemma \ref{lem:boringgens}, we may assume $k_3=0$. We prove this by strong induction on $k_1+k_2$. If $k_1+k_2=0$ or 1, then $k_1=0$ or $k_2=0$. In this case, 1 is not a shared index for the generators of $\mathfrak{a}$ and $\mathfrak{a}'=\mathfrak{a}$. The result follows from Lemmas \ref{lem:oneatatime} and \ref{lem:boringgens}. We do the base case of $k_1=k_2=1$ so that $k_1+k_2=2$. In this case, we have $\mathfrak{a}=(x_1g_1,y_1g_2)$. We write
$\mathfrak{a}=(x_1,y_1) \cap (x_1,g_2) \cap (g_1,y_1g_2).$

Neither piece 2 nor 3 have a minimal prime containing $(x_1,y_1)$, so by Lemma \ref{lem:oneatatime}, we remove the first component and intersect components 2 and 3, leaving us with $(x_1g_1,x_1y_1g_2,[g_1g_2],y_1g_2)$. Since the second generator is divisible by the last, we remove it, leaving us with $(x_1g_1,y_1g_2,[g_1g_2])$.

 For the inductive step, suppose that $k_1+k_2 >2$, $k_1,k_2 \ge 1$, and without loss of generality, that $k_1>1$. We have
\begin{align*}
    \mathfrak{a}&=(x_1,y_1g_{21},\ldots,y_1g_{2k_2}) \cap (g_{11},x_1g_{12},\ldots,x_1g_{1k_1},y_1g_{21},\ldots,y_1g_{2k_2}).
\end{align*}
By Lemma \ref{lem:combinecanonical}, it suffices to recompose these ideals with respect to the index $1$ and intersect the results. By the induction hypothesis, this gives us:
\begin{align*}
    (x_1,g_{21},\ldots,g_{2k_2}) \cap & (g_{11},x_1g_{12},\ldots,x_1g_{1k_1}, y_1g_{21},\ldots,y_1g_{2k_2},\\
    &\qquad [g_{1j_1}g_{2j_2}]_{1 \le j_1 \le k_1, 1 \le j_2 \le k_2}).
\end{align*}
Intersecting these, we get the desired result.
\end{proof}

\begin{rem}
\label{rem:connecttorfstructure}
We can interpret this lemma in terms of the RF-structure of the neural code as well. For each pair of generators $x_1g_1, y_1g_2$, say $g_1=\prod_{i \in \sigma_1} x_i \prod_{i \in \tau_1} y_i$ and $g_2=\prod_{i \in \sigma_2} x_i \prod_{i \in \tau_2} y_i$, where $1 \not\in \sigma_1 \cup \sigma_2 \cup \tau_1 \cup \tau_2$, we have relations
\begin{align*}
    U_1 \bigcap (\bigcap_{i \in \sigma_1} U_i) \subseteq \bigcup_{i \in \tau_1} U_i \quad \text{and} \quad \bigcap_{i \in \sigma_2} U_i \subseteq U_1 \bigcup (\bigcup_{i \in \tau_2} U_i).
\end{align*}
We claim that these two relations together imply the relation
\[(\bigcap_{i \in \sigma_1} U_i) \bigcap (\bigcap_{i \in \sigma_2} U_i ) \subseteq (\bigcup_{i \in \tau_1} U_i) \bigcup (\bigcup_{i \in \tau_2} U_{i}).\]
To see this, we set $X=U_1$, $A=\bigcap_{i \in \sigma_1} U_i$, $B=\bigcup_{i \in \tau_1} U_i$, $C=\bigcap_{i \in \sigma_2} U_i$, and $D=\bigcup_{i \in \tau_2} U_i$. Then $X \cap A \subseteq B$ and $C \subseteq X \cup D$. Since $A \cap C \subseteq C$, it is also contained in $X \cup D$. Then 
\[A \cap C \subseteq A \cap (X \cup D)=(A \cap X) \cup (A \cap D) \subseteq B \cup (A \cap D) \subseteq B \cup D,\]
proving our claim.

As a result, by \cite[Section 4.3]{neuralring} the neural ideal contains a monomial \[[g_1g_2]=\left[\prod_{i \in \sigma_1 \cup \sigma_2} x_i \prod_{i \in \tau_1 \cup \tau_2} y_i \right].\] If $g_1$ and $g_2$ share no indices, then either this monomial or a monomial dividing it must be in the canonical form. If $g_1$ and $g_2$ share an index, this monomial will be removed in Step 5 of Algorithm \ref{alg:canonform}.
\end{rem}

\begin{prop}
\label{pr:origgensstay}
Let $\mathfrak{a}=(g_1,\ldots,g_k) \subseteq S$ be a squarefree monomial ideal. Then \fourstepping\ $\mathfrak{a}$ (with respect to any shared index or to all shared indices) returns an ideal whose generators include $g_1,\ldots,g_k$.

Consequently, for any $j$ such that $x_iy_i \nmid g_j$ for any $1 \le i \le n$, the generators of the \fivestep\ of $\mathfrak{a}$ include $g_j$. Any removal of such a $g_j$ happens in Step 6. Hence for each such $j$, the canonical form of $\mathfrak{a}$ has a generator dividing $g_j$.

In addition, if $\mathfrak{a}$ has a generator $g_j$ such that $x_iy_i \mid g_j$ for some $1 \le i \le n$, $g_j$ may be removed at any stage of Algorithm \ref{alg:canonform}.
\end{prop}

\begin{proof}
To prove that \fourstepping\ $\mathfrak{a}$ does not remove $g_j$ from the ideal, the main thing we need to prove is that removing primes containing an $x_i,y_i$ pair from the primary decomposition of $\mathfrak{a}$ does not remove $g_j$. We prove the result for recomposing with respect to a single index. The result for recomposing with respect to multiple indices follows by induction.

Without loss of generality, suppose that $1$ is a shared index for the generators of $\mathfrak{a}$. Rename the generators of $\mathfrak{a}$ as
\[\mathfrak{a}=(x_1y_1g_{11},\ldots,x_1y_1g_{1k_1},x_1g_{21},\ldots,x_1g_{2k_2},y_1g_{31},\ldots,y_1g_{3k_3},g_{41},\ldots,g_{4k_4}),\]
where $k_j \ge 0$ for $i=1,2,3,4$, and $x_1,y_1 \nmid g_{ij}$ for any $i,j$.
We work by induction on $k_1$.

If $k_1=0$, $\mathfrak{a}=(x_1g_{21},\ldots,x_1g_{2k_2},y_1g_{31},\ldots,y_1g_{3k_3},g_{41},\ldots,g_{4k_4}).$ By Lemma \ref{lem:oneatatime}, recomposing $\mathfrak{a}$ agrees with recomposing the ideal $\mathfrak{a}'$ we get by removing all primes containing both $x_1$ and $y_1$ from $\mathfrak{a}$. By Lemma \ref{lem:lotsofsamexiyi} and Lemma~\ref{lem:boringgens},
\begin{align*}
    \label{eq:oneindexdone}
    \mathfrak{a}'= & (x_1g_{21},\ldots,x_1g_{2k_2}, y_1g_{31},\ldots,y_1g_{3k_3}, \\
    & \quad [g_{2j_2}g_{3j_3}]_{1 \le j_2 \le k_2, 1 \le j_3 \le k_3}, g_{41},\ldots,g_{4k_4}).
\end{align*}
This contains all of the original generators of the ideal. 
Thus when $k_1=0$, the original generators of $\mathfrak{a}$ are generators of the \fourstep\ $\mathfrak{a}'$, and hence of the \fourstep\ $\mathfrak{a}$.

Our induction hypothesis is that for $ \ell_1 = k_1-1$, $\ell_2,\ell_3,\ell_4 \ge 0$, and any choices of $h_{ij}$ not divisible by $x_1$ or $y_1$, recomposing the ideal
\[\mathfrak{a}=(x_1y_1h_{11},\ldots,x_1y_1h_{1\ell},x_1h_{21},\ldots,x_1h_{2\ell_2},y_1h_{31},\ldots,y_1h_{3\ell_3},h_{41},\ldots,h_{4\ell_4})\]
with respect to the index 1 returns
\begin{align*}
   \mathfrak{a}'=(&x_1y_1h_{11},\ldots,x_1y_1h_{1\ell},x_1h_{21},\ldots,x_1h_{2\ell_2},y_1h_{31},\ldots,y_1h_{3\ell_3},\\
   &[h_{2j_2}g_{3j_3}]_{1 \le j_2 \le k_2, 1 \le j_3 \le k_3},h_{41},\ldots,h_{4\ell_4}).
\end{align*}

If $k_1>0$, we have:
\begin{align*}
    \mathfrak{a}=& (x_1,y_1g_{31},\ldots,y_1g_{3k_3},g_{41},\ldots) \cap (y_1,x_1g_{21},\ldots,x_1g_{2k_2},g_{41},\ldots) \\
    &\ \cap (g_{11},x_1y_1g_{12},\ldots,x_1y_1g_{1k_1},x_1g_{21},\ldots,x_1g_{2k_2},y_1g_{31},\ldots,y_1g_{3k_3},g_{41},\ldots).
\end{align*}
By Lemma \ref{lem:combinecanonical}, recomposing $\mathfrak{a}$ with respect to the index 1 agrees with recomposing (with respect to the index 1) and then intersecting the 3 components. Lemma \ref{lem:lotsofsamexiyi} along with the induction hypothesis gives us their recomposed forms:
\begin{align*}
    (x_1,g_{31},\!\ldots\!,g_{3k_3},g_{41},\!\ldots\!)\!\cap\!(y_1,g_{21},\!\ldots\!,g_{2k_2},g_{41},\!\ldots\!)\!\cap\!(g_{11},x_1y_1g_{12},\!\ldots\!,x_1y_1g_{1k_1}, \\ x_1g_{21},\ldots,x_1g_{2k_2},y_1g_{31},\ldots,y_1g_{3k_3}, [g_{2j_2}g_{3j_3}]_{1 \le j_2 \le k_2, 1 \le j_3 \le k_3},g_{41},\ldots,g_{4k_4}).
\end{align*}
Note that in each of the first two pieces, we have already removed a number of pieces divisible by $x_1$ (resp. $y_1$), following Remark \ref{rem:removemultiples}.

Intersecting the first two ideals (while removing monomials divisible by another monomial) gives
\begin{align*}
    (x_1y_1,x_1g_{21},\!\ldots\!,x_1g_{2k_2}, y_1g_{31},\!\ldots\!,y_1g_{3k_3},[g_{2j_2}g_{3j_3}]_{1 \le j_2 \le k_2, 1 \le j_3 \le k_3},g_{41},\!\ldots\!,g_{4k_4}\!).
\end{align*}
Intersecting this with the 3rd piece then gives
\begin{align*}
    (x_1y_1g_{11},\ldots,x_1y_1g_{1k_1},x_1g_{21},\ldots,x_1g_{2k_2},y_1g_{31},[g_{21}g_{31}],\ldots,[g_{2k_2}g_{31}],\\
    \ldots,y_1g_{3k_3},[g_{21}g_{3k_3}],\ldots,[g_{2k_2}g_{3k_3}],g_{41},\ldots,g_{4k_4}),
\end{align*}
as desired.
The original generators of $\mathfrak{a}$ are still generators of this ideal, proving the first statement.

For the second statement, we see by the first statement that every generator $g_j$ of $\mathfrak{a}$ appears in the \fourstep\ $\mathfrak{a}$, and hence if $x_iy_i \nmid g_j$ for any $1 \le i \le n$, $g_j$ appears in the \fivestep\ of $\mathfrak{a}$. As a result, performing Step 6 of Algorithm \ref{alg:canonform} leaves us with a generator that divides $g_j$ for each $j$ such that $x_iy_i \nmid g_j$ for any $1 \le i \le n$.

We return to viewing $\mathfrak{a}=(g_1,\ldots,g_k)$. Now suppose $x_iy_i \mid g_j$ for some $1 \le j \le k$ and $1 \le i \le n$. We see from the proof of the first statement that recomposing $\mathfrak{a}$ with respect to the index $i$ or any other index $i'$ such that $x_{i'}y_{i'} \mid g_j$ preserves $g_j$, but $g_j$ does not contribute to any of the new generators produced by this recomposition. When we perform Step 5 of Algorithm \ref{alg:canonform}, $g_j$ will be removed.

We have 3 cases for what happens to $g_j$ when we recompose $\mathfrak{a}$ with respect to an index $i'$ such that $x_{i'}y_{i'} \nmid g_j$: first, if $x_{i'},y_{i'} \nmid g_j$, preserves $g_j$, and $g_j$ does not contribute to any of the new generators produced by this recomposition. When we perform Step 5 of Algorithm \ref{alg:canonform}, $g_j$ will be removed since it is divisible by $x_iy_i$.

Second, if $x_{i'} \mid g_j$, but $y_{i'} \nmid g_j$, the recomposition produces $g_j$ and some number of $\left[\frac{g_j}{x_{i'}}h \right]$. All of these are removed in Step 5 of Algorithm \ref{alg:canonform} as they are divisible by $x_iy_i$. The third case, where $y_{i'} \mid g_j$ and $x_{i'} \nmid g_j$, is similar.

Since any generator built from $g_j$ will be removed in Step 5 of Algorithm \ref{alg:canonform}, removing $g_j$ before recomposing $\mathfrak{a}$ with respect to any index or set of indices does not change the result of Algorithm \ref{alg:canonform}.
\end{proof}

\begin{rem}
\label{rem:recomposednotunique}
By Remark \ref{rem:removemultiples} and Proposition \ref{pr:origgensstay}, we may remove generators divisible by $x_iy_i$ for some $1 \le i \le n$ or by another generator at any point of Algorithm \ref{alg:canonform} without altering the result. We use this to shorten computations throughout Sections \ref{sec:mainresults}, \ref{sec:genericcanonicalforms}, and \ref{sec:computations}. Consequently, we refer to an ideal that only requires Steps 5 and 6 (resp. Step 6) of Algorithm \ref{alg:canonform} to find the canonical form as recomposed (resp. in \fivestep), even though these ideals may not agree with the results of Step 4 (resp. Step 5)
of the algorithm.
\end{rem}

\section{Classification of canonical neural ideals}
\label{sec:mainresults}

In this section, we apply the results of Section \ref{sec:keylemmas} to determine which polarized neural ideals are in canonical form. Our main result, Theorem \ref{thm:mainthm}, gives a simple criterion for determining whether a neural ideal is in canonical form based only on the shared indices in its generators. Notably, our results do not depend on the total number of neurons, only the indices involved in the generators of the ideal.

\begin{rem}
\label{rem:onegen}
Any neural ideal with a single generator is in canonical form. This is because the components of its primary decomposition are all principal, and cannot contain $(x_i,y_i)$ for any $1 \le i \le n$. Alternatively, with only one generator, there is no pair of generators that can share an index. 
\end{rem}

\begin{thm}
\label{thm:2gens}
The canonical forms of polarized neural ideals with 2 generators, neither dividing the other, are as follows:
\begin{enumerate}
    \item Any polarized neural ideal of the form $(g_1,g_2)$ where $g_1$ and $g_2$ share no indices, is in canonical form.
    \item Any polarized neural ideal of the form \[(x_1\cdots x_tg_1,y_1 \cdots y_tg_2),\] where $t \ge 2$ and $g_1$ and $g_2$ share no indices, is in canonical form.
    \item Any polarized neural ideal of the form \[(x_1 \cdots x_ty_{t+1}\cdots y_{t+s}g_1,x_{t+1}\cdots x_{t+s}y_1\cdots y_tg_2),\] where $t,s>0$ and $g_1$ and $g_2$ share no indices, is in canonical form.
    \item Any polarized neural ideal of the form $(x_1g_1,y_1g_2)$, where $g_1$ and $g_2$ share no indices, is not in canonical form.
    Recomposing the ideal yields $(x_1g_1,y_1g_2,[g_1g_2])$, and so the canonical form of this neural ideal is:
    \begin{enumerate}
        \item $(x_1g_1,y_1g_2,[g_1g_2])$ if $[g_1g_2] \ne g_1,g_2$,
        \item $(x_1g_1,g_2)$ if $[g_1g_2]=g_2 \ne g_1$,
        \item $(g_1,y_1g_2)$ if $[g_1g_2]=g_1 \ne g_2$
        \item and $(g_1)$ if $g_1=g_2$.
    \end{enumerate}
\end{enumerate}
\end{thm}

\begin{proof}
(1) This follows from Lemma \ref{lem:boringgens}.

(2) We prove this by induction on $k$. First assume $k=2$. We use Lemma \ref{lem:oneatatime} to simplify the computation. We have
    \begin{align*}
        (x_1x_2g_1,y_1y_2g_2)= (x_1,y_1) \cap (x_2g_1,y_1) \cap (x_1,y_2g_2) \cap (x_2g_1,y_2g_2).
    \end{align*}
    The first piece is the only piece containing $(x_1,y_1)$, so we remove it.
    Next we look for components containing $(x_2,y_2)$:
    \begin{align*}
        (x_2g_1,y_1) \cap (x_1,y_2g_2) \cap (x_2g_1,y_2g_2) \hspace{-4cm} & \\
        &=(x_2g_1,y_1) \cap (x_1,y_2g_2) \cap (x_2,y_2) \cap (x_2,g_2) \cap (g_1,y_2) \cap (g_1,g_2).
    \end{align*}
    Removing $(x_2,y_2)$ and completing Algorithm \ref{alg:canonform} yields $(x_1x_2g_1,y_1y_2g_2)$, as desired.
    Now suppose $t>2$. We have
    \begin{align*}
        &(x_1 \cdots x_tg_1,y_1 \cdots y_tg_2)=(x_1,y_1) \cap (x_1,y_2 \cdots y_tg_2) \cap (x_2 \cdots x_tg_1,y_1 \cdots y_tg_2).
    \end{align*}
    By Lemma \ref{lem:oneatatime}, the canonical form of the original ideal is equal to the canonical form of the intersection of pieces 2 and 3. By (1) and the induction hypothesis, both pieces are in canonical form. Hence by Lemma \ref{lem:combinecanonical}, we may apply Steps 5 and 6 of Algorithm \ref{alg:canonform} to their intersection to get the canonical form of the original ideal. Intersecting the ideals 
    \[(x_1,y_2 \cdots y_tg_1) \cap (x_2 \cdots x_tg_1,y_1 \cdots y_tg_2)\]
    and completing Steps 5 and 6 of Algorithm \ref{alg:canonform}, we get the desired result.
    
(3) We prove this by induction on $t+s$, starting with the base case $t=s=1$. In this case, we have
    \begin{align*}
        &(x_1y_2g_1,x_2y_1g_2)=(x_1,y_1) \cap (x_1,x_2g_2) \cap (y_2,x_2)\cap (g_1,x_2) \cap (y_2g_1,y_1g_2).
    \end{align*}
    Removing the pieces $(x_1,y_1)$ and $(y_2,x_2)$ as in Lemma \ref{lem:oneatatime}, we intersect the rest of the ideals following Algorithm \ref{alg:canonform} to get the ideal $(x_1y_2g_1,x_2y_1g_2)$.
    
    Now assume that $t+s>1$, and without loss of generality that $t>1$. We have
    \begin{align*}
        (x_1 \cdots x_ty_{t+1}\cdots y_{t+s}g_1,x_{t+1}\cdots x_{t+s}y_1 \cdots y_tg_2) \hspace{-5cm} & \\
        &=(x_1,x_{t+1}\cdots x_{t+s}y_1 \cdots y_tg_2)\cap (x_2 \cdots x_ty_{t+1} \cdots y_{t+s}g_1,x_{t+1}\cdots x_{t+s}y_1 \cdots y_tg_2).
    \end{align*}
    By the induction hypothesis, the second piece is in canonical form (1 is no longer a shared index). By part (4), the canonical form of the first piece is $(x_1,x_{t+1} \cdots x_{t+s}y_2 \cdots y_tg_2)$. Combining these using Lemma \ref{lem:combinecanonical}, we get the desired result.
    
(4) We have $(x_1g_1,y_1g_2)=(x_1,y_1) \cap (x_1,g_2) \cap (g_1,y_1g_2).$
    By Lemma \ref{lem:oneatatime}, the canonical form of the original ideal is equal to the canonical form of the intersection of the second and third ideals.
    Removing the first piece and combining the other two via Algorithm \ref{alg:canonform}, we get the \fivestep\ $(x_1g_1,y_1g_2,[g_1g_2])$. The result follows by removing generators that are a multiple of another generator. If $[g_1g_2] \ne g_1,g_2$, we remove no generators. If $[g_1g_2]=g_1\ne g_2$, we remove $x_1g_1$. If $[g_1g_2]=g_2 \ne g_1$, we remove $y_1g_2$. If $[g_1g_2]=g_1=g_2$, we remove $x_1g_1$ and $y_1g_2$. \qedhere
\end{proof}

\begin{example}
We give examples of each case in Theorem \ref{thm:2gens}.
\begin{enumerate}
    \item The neural ideal $(x_1y_2,x_3y_2)$ is in canonical form.
    \item The neural ideal $(x_1x_2y_3,x_4y_1y_2)$ is in canonical form.
    \item The neural ideal $(x_1y_2,x_2y_1)$ is in canonical form.
    \item The neural ideal $(x_1y_2,x_3y_1)$ has canonical form $(x_1y_2,x_3y_1,x_3y_2)$. In contrast, the neural ideal $(x_1,x_3y_1)$ has canonical form $(x_1,x_3)$.
\end{enumerate}
\end{example}

\begin{rem}
Following the lead of Remark \ref{rem:connecttorfstructure}, in the case where the two generators share a single index and either $g_1=1$ or $g_2=1$,  we can see that the ideal is not in canonical form more directly from the generators and the conditions for the canonical form given in \cite[Theorem 4.3]{neuralring}:

First, we assume that $g_2=1$, so the ideal is $(x_1g_1,y_1)$. Suppose that $g_1 =x_2 \cdots x_ty_{t+1}\cdots y_{t+s}$. Then the generators tell us that 
\begin{equation}
    \label{eq:contradiction1}
    U_{1\cdots t} \subseteq U_{t+1}\cup \ldots \cup U_{t+s}, \, X \subseteq U_1.
\end{equation}
But then $U_{1 \cdots t}=U_{2 \cdots t},$
so the left hand side in Equation \ref{eq:contradiction1} is non-minimal. Hence $(x_1g_1,y_1)$ is not in canonical form. Further, removing $x_1$ from the first generator fixes the issue, which agrees with the result of Theorem \ref{thm:2gens}.

Now suppose that $g_1=1$, so that the ideal becomes $(x_1,y_1g_2)$. Suppose that $g_2=x_2 \cdots x_ty_{t+1}\cdots y_{t+s}$. Then the generators tell us that
\begin{equation}
    \label{eq:contradiction2}
    U_1=\emptyset, \, U_{2 \cdots t} \subseteq U_1 \cup U_{t+1} \cup \ldots \cup U_{t+s}.
\end{equation}
Since $U_1=\emptyset$, the right hand side in Equation \ref{eq:contradiction2} is non-minimal. Hence $(x_1,y_1g_2)$ is not in canonical form. In this case, removing $y_1$ from the second generator fixes the issue, which is consistent with the result of Theorem \ref{thm:2gens}.
\end{rem}

\begin{thm}
\label{thm:mainthm}
Let $\mathfrak{a}=(g_1,\ldots,g_k)$ be a polarized neural ideal such that $g_{j_1} \nmid g_{j_2}$ for any $1 \le j_1 \ne j_2 \le k$, and $x_iy_i \nmid g_j$ for any $1 \le i \le n$ and $1 \le j \le k$. 
If for some pair $g_{j_1},g_{j_2}$ of generators of $\mathfrak{a}$, $g_{j_1}$ and $g_{j_2}$ share exactly 1 index $i$ and no other generator of $\mathfrak{a}$ divides $\frac{[g_{j_1}g_{j_2}]}{x_iy_i}$, then $\mathfrak{a}$ is not in canonical form.
Otherwise, $\mathfrak{a}$ is in canonical form.
\end{thm}

\begin{proof}
First assume $\mathfrak{a}$ has at least one pair of generators with exactly one shared index, say the index 1. To distinguish the generators with a factor of $x_1$ or $y_1$, we rewrite $\mathfrak{a}=(x_1g_{11},\ldots,x_1g_{1k_1},y_1g_{21},\ldots,y_1g_{2k_2},g_{31},\ldots,g_{3k_3}),$
where $x_1,y_1 \nmid g_{j\ell}$ for any $j,\ell$. 
By Lemma \ref{lem:lotsofsamexiyi}, recomposing $\mathfrak{a}$ with respect to the index 1 returns
\begin{align*}
    \mathfrak{a}'=(&x_1g_{11},\!\ldots\!,x_1g_{1k_1},y_1g_{21},\!\ldots\!,y_1g_{2k_2},[g_{1j_1}g_{2j_2}]_{1 \le j_1 \le k_1, 1 \le j_2 \le k_2},g_{31},\!\ldots\!,g_{3k_3}),
\end{align*}
which by Lemma \ref{lem:oneatatime} has the same canonical form of $\mathfrak{a}$.

By our hypothesis, some pair $g_{1j_1}$, $g_{2j_2}$ has no shared indices, so by Proposition \ref{pr:origgensstay}, the canonical form of $\mathfrak{a}$ either has the generator $[g_{1j_1}g_{2j_2}]$ or has another generator $g$ that is a proper divisor of $[g_{1j_1}g_{2j_2}]$.

If no original generator of $\mathfrak{a}$ divides $[g_{1j_1}g_{2j_2}]$, then $g$ is not an original generator of $\mathfrak{a}$, so $\mathfrak{a}$ is not in canonical form.

If, however, for every index $i$ and every pair of generators $g_{j_1}$ and $g_{j_2}$ of $\mathfrak{a}$ that share only the index $i$, some other generator of $\mathfrak{a}$ divides $\frac{[g_{j_1}g_{j_2}]}{x_iy_i}$, then every term of the form $\frac{[g_{j_1}g_{j_2}]}{x_iy_i}$ is removed by Step 6 of Algorithm \ref{alg:canonform}, and by Remark \ref{rem:removemultiples} may be removed now. Hence recomposing $\mathfrak{a}$ returns $\mathfrak{a}$. By the hypotheses that no $x_iy_i$ divides any $g_j$ and no $g_j$ divides another, $\mathfrak{a}$ is in canonical form.

Next assume that no generators of $\mathfrak{a}$ share any indices. Then $\mathfrak{a}$ is in canonical form by Lemma \ref{lem:boringgens}.

Finally, assume that no generators of $\mathfrak{a}$ share exactly one index, and that at least one pair of generators shares an index. 
Without loss of generality, assume 1 is a shared index among the generators of $\mathfrak{a}$. Rewrite 
\[\mathfrak{a}=(x_1g_{11},\ldots,x_1g_{1k_1},y_1g_{21},\ldots,y_1g_{2k_2},g_{31},\ldots,g_{3k_3}),\]
where $x_1,y_1 \nmid g_{j\ell}$ for any $j,\ell$.

By Lemma \ref{lem:lotsofsamexiyi}, recomposing $\mathfrak{a}$ with respect to the index 1 returns
\begin{align*}
    \mathfrak{a}'=(&x_1g_{11},\!\ldots\!,x_1g_{1k_1},y_1g_{21},\!\ldots\!,y_1g_{2k_2},[g_{1j_1}g_{2j_2}]_{1 \le j_1 \le k_1, 1 \le j_2 \le k_2},g_{31},\!\ldots\!,g_{3k_3}),
\end{align*}
which by Lemma \ref{lem:oneatatime} has the same canonical form as $\mathfrak{a}$.

However, since any two generators that share an index share at least 2 indices, each $[g_{1j_1}g_{2j_2}]$ is divisible by $x_iy_i$ for some $2 \le i \le n$. So they will all vanish when Step 5 of Algorithm \ref{alg:canonform} is applied. By Proposition \ref{pr:origgensstay}, we can remove them now, leaving us with $\mathfrak{a}$ again. We repeat this process for each additional shared index, getting $\mathfrak{a}$ again each time. Then $\mathfrak{a}$ recomposed is equal to $\mathfrak{a}$. By assumption, $x_iy_i$ does not divide any generator of $\mathfrak{a}$, so $\mathfrak{a}$ is then in \fivestep. Since we assumed no generator of $\mathfrak{a}$ divides another, $\mathfrak{a}$ is in canonical form.
\end{proof}

Note that the divisibility condition in Theorem \ref{thm:mainthm} requires the ideal to have at least 3 generators, hence why it is not a part of Theorem \ref{thm:2gens}.

\begin{example}
We include an example to illustrate the weirdest case of Theorem \ref{thm:mainthm}. Let $\mathfrak{a}=(x_1x_2,x_3x_4y_1,x_2x_3)$. The first two generators share only the index 1, but this ideal is nevertheless in canonical form. By Theorem \ref{thm:2gens} along with Lemma \ref{lem:boringgens}, the \fivestep\ of this ideal is
\[(x_1x_2,x_3x_4y_1,x_2x_3x_4,x_2x_3).\]
We remove $x_2x_3x_4$ since it is divisible by $x_2x_3$, returning the original ideal. Hence $\mathfrak{a}$ is in canonical form.
\end{example}

\begin{cor}
Let $\mathfrak{a}$ be a polarized neural ideal such that any pair of generators that share the index $i$ share at least 2 indices. Then recomposing $\mathfrak{a}$ with respect to the index $i$ returns $\mathfrak{a}$.

If each pair of generators that share any index share at least 2 indices, then recomposing $\mathfrak{a}$ with respect to any set of indices returns $\mathfrak{a}$.
\end{cor}

\begin{proof}
This follows from the proof of Theorem \ref{thm:mainthm}.
\end{proof}

\begin{cor}
\label{cor:sharedindexrecomposition}
Let $\mathfrak{a}$ be a polarized neural ideal. In order to compute the canonical form of $\mathfrak{a}$, it suffices to recompose $\mathfrak{a}$ with respect to the indices $i_1,\ldots,i_t$ such that for each $1 \le j \le t$, some pair of generators of $\mathfrak{a}$ shares only the index $i_j$.
\end{cor}

\begin{proof}
By Lemma \ref{lem:lotsofsamexiyi}, recomposing 
\[\mathfrak{a}=(x_ig_{11},\ldots,x_ig_{1k_1},y_ig_{21},\ldots,y_ig_{2k_2},g_{31},\ldots,g_{3k_3})\] with respect to the index $i$ returns
\begin{align*}
    \mathfrak{a}'=(x_ig_{11},\ldots,x_ig_{1k_1}, y_ig_{21},\ldots,y_ig_{2k_2}, [g_{1j_1}g_{2j_2}]_{1 \le j_1 \le k_2, 1 \le j_2 \le k_2}, g_{31},\ldots,g_{3k_3}).
\end{align*}
As in the proof of Theorem \ref{thm:mainthm}, if each pair $x_ig_{1j_1},y_ig_{2j_2}$ sharing the index $i$ shares an additional index, every $[g_{1j_1}g_{2j_2}]$ can be removed by Proposition \ref{pr:origgensstay}, and we get $\mathfrak{a}$ back.
\end{proof}

We can now give a new algorithm for computing the canonical form:

\begin{algorithm}
\label{alg:shorter}
Begin with a polarized neural ideal $\mathfrak{a}=(g_1,\ldots,g_k)$.
\begin{enumerate}
    \item List all indices $\{i_1,\ldots,i_{\ell}\}$ such that some pair of generators of $\mathfrak{a}$ share only that index. If there are none, jump to the last step.
    \item Starting with the first index $i_1$ in the list, for each pair of generators $g_{j_1},g_{j_2}$ of $\mathfrak{a}$ that share only the index $i_1$, add a generator $\frac{[g_{j_1}g_{j_2}]}{x_{i_1}y_{i_1}}$ . Replace the generators of $\mathfrak{a}$ with this new longer list of generators.
    \item Repeat for $i_2,\ldots,i_{\ell}$.
    \item When Step 2 has been completed for every index in the original list, remove any remaining generators divisible by $x_iy_i$ for any $1 \le i \le n$ and any generators divisible by another generator of smaller degree.
\end{enumerate}
\end{algorithm}

\begin{proof}
By Corollary \ref{cor:sharedindexrecomposition}, it suffices to recompose $\mathfrak{a}$ with respect to the indices $i$ such that some pair of generators $g_{j_1}$ and $g_{j_2}$ of $\mathfrak{a}$ share only the index $i$. This is the list made in Step 1. Performing Step 2 recomposes $\mathfrak{a}$ with respect to the index $i$ by Lemma \ref{lem:lotsofsamexiyi}, with one change: any new generator coming from a pair of generators $g_{j_1},g_{j_2}$ that share more than one index is skipped, as by the proof of Theorem \ref{thm:mainthm} they would be removed in Step 5 of Algorithm \ref{alg:canonform} and by Proposition \ref{pr:origgensstay} they may be removed at any time.
Hence after Step 3, $\mathfrak{a}$ is recomposed.
 Step 4 removes any remaining generators divisible by $x_iy_i$ or by another generator that have not been removed already, completing Steps 5 and 6 of Algorithm \ref{alg:canonform}.
\end{proof}

\section{Generic canonical forms}
\label{sec:genericcanonicalforms}

The results of this section give us a way to compare canonical forms of neural ideals where the components with no shared indices have been changed. We do this by passing to an extension ring, computing a ``generic" canonical form, and then returning to our original ring.

\begin{notation}
Let $g_1,\ldots,g_k \in S$ be squarefree monomials sharing no indices. We define an extension ring $P=S[z_1,\ldots,z_k]$ and a map $\pi:P \to S$ sending $z_j \mapsto g_j$. Then $S \cong P/(z_j-g_j)$.
\end{notation}

We compute an analogue of the canonical form for squarefree monomial ideals in $P$ by modifying Algorithm \ref{alg:canonform} as follows:

\begin{algorithm}
\label{alg:polarizemore}
Let $P=S[z_1,\ldots,z_k]$.
\begin{enumerate}
    \item Start with a squarefree monomial ideal $\mathfrak{a}=(f_1,\ldots,f_m)$ in $P$, where $f_t=\prod_{i \in \sigma} x_i \prod_{i \in \tau} y_i \prod_{i \in \upsilon} z_i$.
    \item Compute the primary decomposition of $\mathfrak{a}$. The ideals $p_1,\ldots,p_s$ in the primary decomposition will all be generated by a subset of $(x_1,\ldots,x_n,y_1,\ldots,y_n,z_1,\ldots,z_k)$.
    \item Impose the relation $x_i+y_i=1$, and as a result remove any ideal in the primary decomposition containing both $x_i$ and $y_i$ for any $1 \le i \le n$.
    \item Compute $[h_1\cdots h_s]$ for every set of choices of $h_{\ell}$ a generator of $p_{\ell}$. (That is, intersect the remaining prime ideals.)
    \item Impose the relations $x_iy_i=0$, and remove any monomials with a factor of $x_iy_i$ for some $1 \le i \le n$.
    \item Remove any remaining products that are multiples of a product of lower degree.
\end{enumerate}
\end{algorithm}

\begin{defn}
We say that an ideal that is the result of Algorithm \ref{alg:polarizemore} is in \textit{generic canonical form}. As we did in $S$, we will refer to an ideal in $P$ that has had Steps 1-4 of Algorithm \ref{alg:polarizemore} applied to it as \textit{generically recomposed}, and an ideal that has had Steps 1-5 of Algorithm \ref{alg:polarizemore} as being in \textit{generic almost canonical form}.
\end{defn}

As with their non-generic version (see Remark \ref{rem:recomposednotunique}), we refer to an ideal that only requires Steps 5 and 6 (resp. Step 6) of Algorithm \ref{alg:polarizemore} to find the generic canonical form as generically recomposed (resp. in generic almost canonical form), even though these ideals may not agree with the results of Step 4 (resp. Step 5) of Algorithm \ref{alg:polarizemore}.

\begin{lemma}
\label{lem:gsdontmatter}
Let $b_1,\ldots,b_k,g_1,\ldots,g_k \in S$ be squarefree monomials such that the $g_j$ share no indices with the $b_jg_j$. 
In recomposing $\mathfrak{a}=(b_1g_1,\ldots,b_kg_k)$, it is not necessary to split the $g_j$ into their factors.
\end{lemma}

\begin{proof}
We proceed by induction on $k$. We noted in Remark \ref{rem:onegen} that neural ideals with one generator are in canonical form, so we begin with $k=2$. In this case we have $(b_1g_1,b_2g_2)=(b_1,b_2) \cap (b_1,g_2) \cap (g_1,b_2) \cap (g_1,g_2).$
By Lemma \ref{lem:boringgens}, since no $g_j$ shares any index with any $b_jg_j$, the second, third, and fourth pieces are already recomposed, and by Lemma \ref{lem:combinecanonical}, recomposed $\mathfrak{a}$ is equal to the intersection of the recomposed components on the right hand side. This proves the case $k=2$.

Now assume $k>2$. We have
\[(b_1g_1,\ldots,b_kg_k)=(b_1,b_2g_2,\ldots,b_kg_k) \cap (g_1,b_2g_2,\ldots,b_kg_k).\]
By Lemma \ref{lem:boringgens}, the recomposed second piece is equal to $g_1$ plus the recomposed $(b_2g_2,\ldots,b_kg_k)$. By the induction hypothesis, we do not need to split the $g_j$ into their factors to compute this. We continue to split the first piece:
\[(b_1,b_2g_2,\ldots,b_kg_k)=(b_1,b_2,b_3g_3,\ldots,b_kg_k) \cap (b_1,g_2,b_3g_3,\ldots,b_kg_k).\]
As before, we do not need to split the $g_j$ into their factors to recompose the second piece. Continuing in this way, we eventually get
\[(b_1,b_2,\ldots,b_{k-1},b_kg_k)=(b_1,\ldots,b_{k-1},b_k) \cap (b_1,\ldots,b_{k-1},g_k).\]
The first piece has no $g_k$, and as before, we do not need to split $g_k$ into its factors to recompose the second piece.
\end{proof}

\begin{notation}
Let $(b_1z_1,\ldots,b_kz_k)$ be a squarefree monomial ideal of $P$. When referring to the generically recomposed, almost canonical, or canonical form of $(b_1z_1,\ldots,b_kz_k)$, we will denote it \[(f(z_1,\ldots,z_k)) := (f_1(z_1,\ldots,z_k), \ldots f_{\ell}(z_1,\ldots,z_k))\] for some $\ell$ where each generator $f_j(z_1,\ldots,z_k) \in P$ is a monomial. This allows us to view the generically recomposed (resp. almost canonical or canonical) form as a function of the $z_j$ that is determined by the $b_j$.
\end{notation}

\begin{prop}
\label{pr:gsdontmatter}
Let $b_1,\ldots,b_k,g_1,\ldots,g_k \in S$ be squarefree monomials such that the $g_j$ share no indices with the $b_jg_j$. Let $P=S[z_1,\ldots,z_k]$ and $\pi:P \to S$ be the homomorphism sending $S$ to itself and $z_j \mapsto g_j$. If the generically recomposed form of $(b_1 z_1, \ldots, b_k z_k)$ in $P$ is given by $(f(z_1,\ldots,z_k))$, then the recomposed form of $(b_1g_1, \ldots, b_k g_k)$ in $S$ is given by \[[\pi(f(z_1,\ldots,z_k))] := ([f_1(g_1,\ldots, g_k)], \ldots, [f_{\ell}(g_1,\ldots,g_k)]).\]

Further, if $x_i y_i$ divides $[f_{j}(g_1,\ldots,g_k)]$ for some $1 \le i \le n$ and $1 \le j \le \ell$, then $x_i y_i$ divides $f_j(z_1,\ldots,z_k)$. 
Hence the canonical form of $(b_1g_1,\ldots,b_kg_k)$ in $S$ may be computed by computing the generic almost canonical form of $(b_1z_1,\ldots,b_kz_k)$, applying $\pi$ and replacing any $x_i^c$ for $c >1$ with $x_i$, and then applying Step 6 of Algorithm \ref{alg:canonform}.

In particular, if we know the canonical form of $(b_1g_1,\ldots,b_kg_k)$ (where the $g_j$ still share no indices with the $b_jg_j$), and $g \mid g_j$, we can find the canonical form of $(b_1g_1,\ldots,b_{j-1}g_{j-1},b_jg,b_{j+1}g_{j+1},\ldots,b_kg_k)$ by computing the canonical form of $(b_1g_1,\ldots,b_kg_k)$, replacing $g_j$ by $g$, and then applying Step 6 of Algorithm \ref{alg:canonform} again.
\end{prop}

\begin{proof}
Say that the generically recomposed $(b_1z_1,\ldots,b_kz_k)$ in $P$ is given by 
\[(f(z_1,\ldots,z_k)):=(f_1(z_1,\ldots,z_k),\ldots,f_\ell(z_1,\ldots,z_k)),\]
where each $f_s(z_1,\ldots,z_k) \in P$.
Applying $\pi$, we get \[(f_1(g_1,\ldots,g_k),\ldots,f_\ell(g_1,\ldots,g_k)).\] We claim that after replacing products with LCMs this is equal to the recomposed $(b_1g_1,\ldots,b_kg_k)$.
By Lemma \ref{lem:gsdontmatter}, since the $g_j$ contain no $x_i$ or $y_i$ with shared indices, we can leave the $g_j$ alone when decomposing $(b_1g_1,\ldots,b_kg_k)$.  Hence $\pi$ applied to the primary decomposition of $(b_1z_1,\ldots,b_kz_k)$ gives us a sufficient decomposition of $(b_1g_1,\ldots,b_kg_k)$ to recompose the latter. Since $\pi(x_i)=x_i$ and $\pi(y_i)=y_i$, the components we remove in Step 3 of Algorithm \ref{alg:canonform} in $S$ are equal to $\pi$ applied to the components we remove in Step 3 of Algorithm \ref{alg:polarizemore}. The only potential difference is in Step 4, and this is resolved by replacing products with LCM's.

Since the $g_j$ contain no indices shared with the $b_jg_j$ (i.e., if $x_i \mid g_\ell$ for any $\ell$, then $y_i \nmid \Pi_j b_jg_j$, and the same with $x_i$ and $y_i$ reversed), this will hold after recomposing $(b_1g_1,\ldots,b_kg_k)$ as well. As a result, any generator of the recomposed $(b_1g_1,\ldots,b_kg_k)$ divisible by $x_iy_i$ for some $1 \le i \le n$ must be equal to $\pi$ of a generator of the recomposed $(b_1z_1,\ldots,b_kz_k)$ divisible by the same $x_iy_i$, up to replacing products with LCMs.
Applying Step 6 of Algorithm \ref{alg:canonform} to $(f_1(g_1,\ldots,g_k),\ldots,f_\ell(g_1,\ldots,g_k))$ now gives us the canonical form of $(b_1g_1,\ldots,b_kg_k)$, as desired.

For the final statement, note that 
\begin{align*}
    \mathfrak{a} &=(b_1g_1,\ldots,b_1g_k)  \text{ and } 
    \mathfrak{a}'=(b_1g_1,\ldots,b_jg,\ldots,b_kg_k)
\end{align*} 
have the same generic canonical form, $(f(z_1,\ldots,z_k))$. Applying $\pi_1$ sending $z_j \mapsto g_j$ and $\pi_2$ sending $z_j \mapsto g$, we have effectively replaced $g_j$ by $g$ everywhere it appears. The only possible difference is in Step 6 of Algorithm \ref{alg:canonform}.
\end{proof}

\begin{example}
\label{eg:genericcanon}
Consider $\mathfrak{a}=(x_1x_2x_4,x_3x_4y_1)$. Set $z_1=\im(x_2x_4)$ and $z_2=\im(x_3x_4)$.
Write $\tilde{\mathfrak{a}}=(x_1z_1,y_1z_2) \subseteq P$. Using the proof of Theorem~\ref{thm:2gens}, this ideal has canonical form $(x_1z_1,y_1z_2,z_1z_2)$ in $P$. Passing back to $S$,
we get $(x_1x_2x_4,x_3x_4y_1,x_2x_3x_4)$. In this case, no generator is a multiple of another, so we are done. This computation agrees with the result of Theorem~\ref{thm:2gens}.

Further, by the last statement of the Theorem, the canonical form of $\mathfrak{b}=(x_1x_2,x_3x_4y_1)$ is $(x_1x_2,x_3x_4y_1,x_2x_3x_4)$ and the canonical form of $\mathfrak{c}=(x_1,x_3x_4y_1)$ is $(x_1x_2,x_3x_4)$. All of these agree with the result of Theorem~\ref{thm:2gens}.
\end{example}

\begin{thm}
\label{thm:repeats}
Suppose that the generic \fivestep\ of the ideal $(b_1z_1,\ldots,b_kz_k)$ in $P=S[z_1,\ldots,z_k]$ is 
\[(f(z_1,\ldots,z_k)):=(f_1(z_1,\ldots,z_k),\ldots,f_{\ell}(z_1,\ldots,z_k)).\]
Let $\mathfrak{a}=(b_1g_{11},\ldots,b_1g_{1t_1},b_2g_{21},\ldots,b_2g_{2t_2},\ldots,b_kg_{k1},\ldots,b_kg_{kt_k}),$ where for any index $i$ shared by the generators of $\mathfrak{a}$, $x_i,y_i \nmid \prod g_{j_1j_2}$.
Then the \fivestep\ of $\mathfrak{a}$ is
\begin{align*}
    (&[g_{1j_1}g_{2j_2}\cdots g_{kj_k}f(1,\ldots,1)],\\
    &[g_{2j_2}\cdots g_{kj_k}f(g_{1t_1},1,\ldots,1)],\ldots,[g_{1j_1}\cdots g_{(k-1)j_{k-1}}f(1,\ldots,1,g_{kt_k})],\\
    &\vdots \\
    &[g_{kj_k}f(g_{1t_1},\ldots,g_{(k-1)t_{k-1}},1)],\ldots,[g_{1j_1}f(1,g_{2t_2},\ldots,g_{kt_k})],\\
    &[f(g_{1t_1},g_{2t_2},\ldots,g_{kt_k})]),
\end{align*}
where $j_i$ ranges over the integers $1,2,\ldots,t_i-1$ and all products have been replaced with LCM's.
\end{thm}

\begin{rem}
The statement of Theorem \ref{thm:repeats} is asymmetrical in the $g_{ij_i}$, with only the $g_{it_i}$ acting as input to the $f$'s. This is a result of the particular decomposition used in the proof, and the result works equally well if we fix for each $1 \le i \le k$ a choice of $1 \le s_i \le t_i$, replace every $g_{it_i}$ with $g_{is_i}$, and let $j_i$ range over the integers $1,2,\ldots,\widehat{s_i},\ldots,t_i$.
\end{rem}

\begin{proof}
We prove this by induction on $\sum_{i=1}^k t_i$. If $t_1=\ldots=t_k=1$, we are done by Proposition \ref{pr:gsdontmatter}. Now suppose $t_i >1$ for at least one $1 \le i \le k$. Without loss of generality, assume $t_1>1$. We have
\begin{align*}
    (b_1g_{11},\!\ldots\!,b_1g_{1t_1},b_2g_{21},\!\ldots\!)=(b_1,b_2g_{21},\!\ldots\!) \cap (g_{11},b_1g_{12},\!\ldots\!,b_1g_{1t_1},b_2g_{21},\!\ldots\!).
\end{align*}
By Lemma \ref{lem:combinecanonical} and Proposition \ref{pr:origgensstay}, the \fivestep\ of $\mathfrak{a}$ is equal to the result of applying Step 5 of Algorithm \ref{alg:canonform} to the intersection of the almost canonical forms of the two components.
By the induction hypothesis and Proposition \ref{pr:gsdontmatter}, the \fivestep\ of the first piece is
\begin{align*}
    (&[g_{2j_2}\cdots g_{kj_k}f(1,\ldots,1)],\\
    &[g_{3j_3}\cdots g_{kj_k}f(1,g_{2t_2},1,\ldots,1)],\ldots,[g_{2j_2}\cdots g_{(k-1)j_{k-1}}f(1,\ldots,1,g_{kt_k})],\\
    &\vdots \\
    &[g_{kj_k}f(1,g_{2t_2},\ldots,g_{(k-1)t_{k-1}},1)],\ldots,[g_{2j_2}f(1,1,g_{3t_3},\ldots,g_{kt_k})],\\
    &[f(1,g_{2t_2},\ldots,g_{kt_k})]),
\end{align*}
where $j_i$ ranges over the integers $1,2,\ldots,t_i$ for $2 \le i \le k$.
By the induction hypothesis and Lemma \ref{lem:boringgens}, the \fivestep\ of the second piece is
\begin{align*}
    (&g_{11},[g_{1j_1}g_{2j_2}\cdots g_{kj_k}f(1,\ldots,1)],\\
    &[g_{2j_2}\cdots g_{kj_k}f(g_{1t_1},1,\ldots,1)],\ldots,[g_{1j_1}\cdots g_{(k-1)j_{k-1}}f(1,\ldots,1,g_{kt_k})],\\
    &\vdots \\
    &[g_{kj_k}f(g_{1t_1},\ldots,g_{(k-1)t_{k-1}},1)],\ldots,[g_{1j_1}f(1,g_{2t_2},\ldots,g_{kt_k})],\\
    &[f(g_{1t_1},g_{2t_2},\ldots,g_{kt_k})]),
\end{align*}
where $j_i$ ranges over the integers $1,2,\ldots,t_i-1$ for $2 \le i \le k$ and $j_1$ ranges over the integers $2,\ldots,t_1-1$.

We claim that
\begin{align*}
    [[f(h_1,\ldots,h_{i-1},1,h_{i+1},\ldots,h_k)][f(h_1,\ldots,h_{i-1},h_i,h_{i+1},\ldots,h_k)]] \hspace{-5cm} & \\
    &=[f(h_1,\ldots,h_{i-1},h_i,h_{i+1},\ldots,h_k)].
\end{align*}
To see this, note that for any $s$ such that $f_s(z_1,\ldots,z_k)$ has no factor of $z_i$, the two components agree, and for any $s$ such that $f_s(z_1,\ldots,z_k)$ has a factor of $z_i$, say $f_s(z_1,\ldots,z_k)=z_ih$, the first component gives $h$ and the second component gives $h_ih$, so their lcm is $h_ih=f_s(h_1,\ldots,h_{i-1},h_i,h_{i+1},\ldots,h_k)$.

Using this claim to intersect the two pieces, the proof is complete.
\end{proof}

\begin{rem}
In computing the canonical form of $\mathfrak{a}$, if $f_s(z_1,\ldots,z_k)$ contains factors of, without loss of generality, $z_1,\ldots,z_p$ but not $z_{p+1},\ldots,z_k$, then 
\begin{align*}
    f_s(g_{1t_1},\ldots,g_{pt_p},1,\ldots,1) = f_s(g_{1t_1},\ldots,g_{pt_p},g_{(p+1)t_{p+1}},1,\ldots,1)\\
    =\ldots=f_s(g_{1t_1},\ldots,g_{kt_k}).
\end{align*}
(See the proof of Theorem \ref{thm:repeats} for justification of this fact.)
So any term with a factor of $f_s(g_{1t_1},\ldots,g_{pg_p},1,\ldots,1)$ or any of the intermediate lines is a multiple of $f_s(g_{1t_1},\ldots,g_{kt_k})$ and will be removed in the canonical form.
\end{rem}

\begin{example}
We compute the \fivestep\ of \[(x_1g_{11},x_1g_{12},y_1g_{21},y_1g_{22},y_1g_{23}),\] where the $g_{j_1j_2}$ contain no shared indices. Continuing from Example \ref{eg:genericcanon}, the generic canonical form of $(x_1g_1,y_1g_2)$ is $(x_1z_1,y_1z_2,z_1z_2)$. Here
\begin{align*}
    f_1(z_1,z_2)&=x_1z_1, & f_2(z_1,z_2)&=y_1z_2, & f_3(z_1,z_2)&=z_1z_2.
\end{align*} In this example, $j_1=1$ and $j_2 \in \{1,2\}$. Applying Theorem \ref{thm:repeats}, we have
\begin{align*}
    g_{1j_1}g_{2j_2}f_1(1,1)=g_{1j_1}g_{2j_2}x_1, \qquad 
    g_{1j_1}g_{2j_2}f_2(1,1&=g_{1j_1}g_{2j_2}y_1, \\
     g_{1j_1}g_{2j_2}f_3(1,1)=g_{1j_1}g_{2j_2},
\end{align*}
so we will only need generators of the third type, which are $[g_{11}g_{21}]$ and $[g_{11}g_{22}]$.
Similarly,
\begin{align*}
    g_{2j_2}f_1(g_{12},1)&=g_{2j_2}x_1g_{12} &
    g_{2j_2}f_2(g_{12},1)&=g_{2j_2}y_1 &
    g_{2j_2}f_3(g_{12},1)&=g_{2j_2}g_{12},
\end{align*}
from which we get generators $[g_{12}g_{21}],[g_{12}g_{22}],y_1g_{21},$ and $y_1g_{22}$,
\begin{align*}
    g_{1j_1}f_1(1,g_{23})&=g_{1j_1}x_1 &
    g_{1j_1}f_2(1,g_{23})&=g_{1j_1}y_1g_{23} &
    g_{1j_1}f_3(1,g_{23})&=g_{1j_1}g_{23},
\end{align*}
from which we get generators $x_1g_{11}$ and $[g_{11}g_{23}]$,
and
\begin{align*}
    f_1(g_{12},g_{23})&=x_1g_{12} &
    f_2(g_{12},g_{23})&=y_1g_{23} &
    f_3(g_{12},g_{23})&=g_{12}g_{23},
\end{align*}
from which we get generators $x_1g_{12},y_1g_{23},$ and $[g_{12}g_{23}]$.
Hence by Theorem \ref{thm:repeats}, the \fivestep\ of 
$\mathfrak{a}=(x_1g_{11},x_1g_{12},y_1g_{21},y_1g_{22},y_1g_{23})$
is
\begin{align*}
    (&[g_{11}g_{21}],[g_{11}g_{22}],[g_{11}g_{23}],[g_{12}g_{21}],[g_{12}g_{22}],[g_{11}g_{23}],\\
    &x_1g_{11},x_1g_{12},y_1g_{21},y_1g_{22},y_1g_{23}).
\end{align*}
This matches the result of Lemma \ref{lem:lotsofsamexiyi}. In this situation of only 1 shared index, using Lemma \ref{lem:lotsofsamexiyi} might be easier, but as the number of shared indices increases, Theorem \ref{thm:repeats} becomes more useful.
\end{example}

\section{Canonical forms of classes of neural ideals}
\label{sec:computations}

In this section we compute the canonical forms of certain classes of polarized neural ideals. As in Section \ref{sec:mainresults}, the results below do not depend on the number of neurons.

\begin{thm}
\label{thm:chain}
The polarized neural ideal
\[\mathfrak{a}=(x_{1}g_1,x_{2}y_{1}g_2,x_{3}y_{2}g_3,\ldots,x_{{k-1}}y_{{k-2}}g_{k-1},y_{{k-1}}g_k),\] where $k > 2$ and the only shared indices are the $x_{i}$ and $y_{i}$ shown, is not in canonical form. The \fivestep\ of $\mathfrak{a}$ is
\begin{align*}
    (&x_{1}g_1,x_{2}g_1g_2,x_{3}g_1g_2g_3,\ldots,x_{{k-1}}g_1 \cdots g_{k-1},g_1 \cdots g_k, \\
    & \quad x_{2}y_{1}g_2,x_{3}y_{1}g_2g_3,\ldots,x_{{k-1}}y_{1}g_2\cdots g_{k-1},y_{1}g_2\cdots g_k,\\
    & \quad x_{3}y_{2}g_3,\ldots,x_{{k-1}}y_{2}g_3 \cdots g_{k-1},y_{2}g_3 \cdots g_k,\\
    & \quad \vdots \\
    & \quad x_{{k-1}}y_{{k-2}}g_{k-1},y_{{k-2}}g_{k-1}g_k,\\
    & \quad y_{{k-1}}g_k).
\end{align*}

If for any $1 \le t \le t' < k$, $[g_t \cdots g_k] \ne [g_t \cdots g_{t'}]$, and for any $1 < s' \le s \le k$, $[g_1 \cdots g_s] \ne [g_{s'} \cdots g_s]$, then this
is the canonical form of the ideal.

If for some $1 \le s' \le s \le t \le t' \le k$, $[g_1 \cdots g_s]=[g_{s'}\cdots g_s]$ and $[g_t \cdots g_k]=[g_t \cdots g_{t'}]$, where $s$ is minimal, $s'$ is maximal, $t$ is maximal, and $t'$ is minimal with respect to these properties, then the canonical form of $\mathfrak{a}$ is 
\begin{align*}
    (&x_{1}g_1,x_{2}[g_1g_2],\ldots,x_{s-1}[g_1 \cdots g_{s-1}],x_s[g_{s'} \cdots g_s],\ldots,x_{t'-1}[g_{s'} \cdots g_{t'-1}],[g_{s'} \cdots g_{t'}], \\
    & \quad x_{2}y_{1}g_2,x_{3}y_{1}[g_2g_3],\ldots,x_{s-1}y_1[g_2 \cdots g_{s-1}],\\
    & \quad \vdots \\
    & \quad x_{s'}y_{s'-1}g_{s'},\ldots,x_{s-1}y_{s'-1}[g_{s'}\cdots g_{s-1}],\\
    & \quad x_{s'+1}y_{s'}g_{s'+1},x_{s'+2}y_{s'}[g_{s'+1}g_{s'+2}],\ldots,x_{t'-1}y_{s'}[g_{s'+1} \cdots g_{t'-1}],y_{s'}[g_{s'+1} \cdots g_{t'}],\\
    & \quad \vdots\\
    & \quad x_ty_{t-1}g_t,\ldots,x_{t'-1}y_{t-1}[g_t \cdots g_{t'-1}],y_{t-1}[g_t \cdots g_{t'}],\\
    & \quad x_{t+1}y_tg_{t+1},\ldots,x_{k-1}y_t[g_{t+1}\cdots g_{k-1}],y_t[g_{t+1}\cdots g_k],\\
    & \quad \vdots\\
    & \quad x_{{k-1}}y_{{k-2}}g_{k-1},y_{{k-2}}[g_{k-1}g_k],\\
    & \quad y_{{k-1}}g_k).
\end{align*}
\end{thm}

\begin{rem}
If $k=2$, we are in the case of Theorem \ref{thm:2gens}, so we skip that case here. If $k=1$, there are no shared indices and so the ideal is in canonical form by Lemma \ref{lem:boringgens}.

Notice that the canonical form of $\mathfrak{a}$ usually has more generators than $\mathfrak{a}$, often many more. The leftmost column of the canonical form above contains the original generators of the ideal, and all other generators are additional.
\end{rem}

\begin{proof}
We prove this by induction, beginning with the case $k=3$.
We decompose the ideal $(x_1g_1,x_2y_1g_2,y_2g_3)$, removing any component that contains a pair $(x_i,y_i)$ as generators as we go by Lemma \ref{lem:oneatatime}:
\begin{align*}
    &(x_1g_1,x_2y_1g_2,y_2g_3)=(x_1,x_2y_1g_2,y_2g_3) \cap (g_1,x_2y_1g_2,y_2g_3) \\
    & =(x_1,y_1,y_2g_3) \cap (x_1,x_2g_2,y_2g_3) \cap (g_1,x_2,y_2g_3) \cap (g_1,y_1g_2,y_2g_3) \\
    & \quad \to (x_1,x_2,y_2g_3) \cap (x_1,g_2,y_2g_3) \cap (g_1,x_2,y_2) \cap (g_1,x_2,g_3) \cap (g_1,y_1g_2,y_2g_3) \\
    & \quad \to (x_1,x_2,y_2) \cap (x_1,x_2,g_3) \cap (x_1,g_2,y_2g_3) \cap (g_1,x_2,g_3) \cap (g_1,y_1g_2,y_2g_3) \\
    & \quad \to (x_1,x_2,g_3) \cap (x_1,g_2,y_2g_3) \cap (g_1,x_2,g_3) \cap (g_1,y_1g_2,y_2g_3).
\end{align*}
By Lemma \ref{lem:boringgens}, since none of the remaining pieces have any shared indices, each piece is recomposed. By Lemma \ref{lem:combinecanonical}, to get the recomposed $\mathfrak{a}$, it suffices to intersect these ideals.
Intersecting these ideals and removing components divisible by $x_iy_i$ or by another component as in Proposition \ref{pr:origgensstay}, we are left with the \fivestep
$(x_1g_1,x_2g_1g_2,g_1g_2g_3,x_2y_1g_2,y_1g_2g_3,y_2g_3).$

Now assume $k > 3$. We have
\begin{align*}
    &(x_1g_1,x_2y_1g_2,\ldots,x_{k-1}y_{k-2}g_{k-1},y_{k-1}g_k) \\
    & \quad =(x_1,y_1,x_3y_2g_3,\ldots,x_{k-1}y_{k-2}g_{k-1},y_{k-1}g_k) \\
    &\qquad \cap (x_1,x_2g_2,x_3y_2g_3,\ldots,y_{k-1}g_k) \cap (g_1,x_2y_1g_2,x_3y_2g_3,\ldots,y_{k-1}g_k).
\end{align*}
We remove the first piece, which contains the pair $(x_1,y_1)$. By Lemma \ref{lem:combinecanonical}, the recomposed $\mathfrak{a}$ is equal to the intersection of the recomposed second piece with the recomposed third piece. By the induction hypothesis, the generically recomposed second piece is 
\begin{align*}
    (&x_1,x_2g_2,x_3g_2g_3,\ldots,x_{k-1}g_2\cdots g_{k-1},g_2\cdots g_k,\\
    & \quad x_3y_2g_3,\ldots,x_{k-1}y_2 g_3 \cdots g_{k-1},y_2 g_3 \cdots g_k,\\
    & \quad \vdots \\
    & \quad x_{k-1}y_{k-2}g_{k-1},y_{k-2} g_{k-1}g_k,\\
    & \quad y_{k-1} g_k).
\end{align*}
Viewing $y_1g_2$ as a piece sharing no indices with other generators, the recomposed 3rd piece is 
\begin{align*}
    (& g_1,x_2y_1 g_2,x_3y_1 g_2 g_3,\ldots,x_{k-1}y_1 g_2\cdots g_{k-1},y_1 g_2\cdots g_k,\\
    & \quad x_3y_2 g_3,\ldots,x_{k-1}y_2 g_3 \cdots g_{k-1},y_2 g_3 \cdots g_k,\\
    & \quad \vdots \\
    & \quad x_{k-1}y_{k-2} g_{k-1},y_{k-2} g_{k-1} g_k,\\
    & \quad y_{k-1}g_k).
\end{align*}
Intersecting the two pieces and removing generators that are divisible by $x_iy_i$, we get
\begin{align*}
    (& x_1 g_1,x_2 g_1g_2,\ldots,x_{k-1} g_1 \cdots g_{k-1}, g_1 \cdots g_k,\\
    & \quad x_2y_1 g_2,x_3y_1 g_2 g_3,\ldots,x_{k-1}y_1 g_2 \cdots g_{k-1},y_1 g_2 \cdots g_k,\\
    & \quad x_3y_2 g_3,\ldots,x_{k-1}y_2 g_3 \cdots g_{k-1}g,y_2 g_3 \cdots g_k,\\
    & \quad \vdots \\
    & \quad x_{k-1}y_{k-2} g_{k-1},y_{k-2} g_{k-1} g_k,\\
    & \quad y_{k-1} g_k).
\end{align*}
Since by assumption, the $g_j$ share no indices, this is the \fivestep\ of $\mathfrak{a}$.

If $[g_1\cdots g_s] \ne [g_{s'} \cdots g_s]$ for any $1 < s' \le s \le k$ and $[g_t \cdots g_k] \ne [g_t \cdots g_{t'}]$ for any $1 \le t \le t' < k$, it is not possible for any generator to be a multiple of another generator, so we have found the canonical form. Otherwise, we remove generators that are a multiple of another generator to get the result.

Regardless of the value of $[g_1 \cdots g_k]$, it will always be a generator of the ideal. Since $[g_1 \cdots g_k]$ has no factors of $x_i$ or $y_i$ for $1 \le i \le k-1$, whereas every original generator of $\mathfrak{a}$ has a factor of $x_i$ or of $y_i$, the canonical form is distinct from the original ideal.
\end{proof}

\begin{example}
Consider the polarized neural ideal $\mathfrak{a}=(x_1g_1,x_2y_1g_2,y_2g_3)$, where the $g_j$ contain no shared indices. By Theorem \ref{thm:chain},
the almost canonical form of this ideal is
$(x_1g_1,x_2[g_1g_2],[g_1g_2g_3],x_2y_1g_2,y_1[g_2g_3],y_2g_3)$.
Set $g_1=~x_3$, $g_2=x_4$, and $g_3=x_3x_4$, so that $[g_1g_2g_3]=g_3$. Substituting our values of $g_1,g_2,$ and $g_3$ and removing generators that are a multiple of another generator, we are left with 
$(x_1x_3,x_2y_1x_4,x_3x_4).$
This is a degenerate case like in parts of Theorem \ref{thm:2gens}, (4)--instead of adding generators to get the canonical form we change a generator.

If instead we set $g_1=g_2=g_3=x_3$, we find every generator is divisible by $x_3$ thereby making the almost canonical form $(x_3)$. This is similar to the last case of Theorem \ref{thm:2gens}.
\end{example}

\begin{thm}
\label{thm:cycle}
The polarized neural ideal
\[(x_1y_kg_1,x_2y_1g_2,\ldots,x_{k-1}y_{k-2}g_{k-1},x_ky_{k-1}g_k),\]
where $k>2$ and the only shared indices are the $x_i$ and $y_i$ shown, 
is not in canonical form. The canonical form of this ideal is
\begin{align*}
    (&x_1y_kg_1,x_2y_k[g_1g_2],\ldots,x_{k-1}y_k[g_1 \cdots g_{k-1}],\\
    &x_2y_1g_2,x_3y_1[g_2g_3],\ldots,x_ky_1[g_1 \cdots g_k],\\
    &x_3y_2g_3,\ldots,x_ky_2[g_3 \cdots g_k],x_1y_2[g_3 \cdots g_kg_1],\\
    &\vdots \\
    &x_ky_{k-1}g_k,x_1y_{k-1}[g_kg_1],\ldots,x_{k-2}y_{k-1}[g_kg_1 \cdots g_{k-2}]).
\end{align*}

In particular, the canonical form of the polarized neural ideal
\[(x_1y_k,x_2y_1,\ldots,x_ky_{k-1})\]
is $(x_iy_j)_{1 \le i,j \le k, i \ne j}$.
\end{thm}

\begin{rem}
Notice that in contrast to Theorem \ref{thm:chain}, the canonical form is the same regardless of $[g_1 \cdots g_k]$. No choice of values for the $g_j$ can force one generator to be a multiple of another, as any two generators have distinct factors of $x_{i_1}y_{i_2}$.

We assume $k>2$ as when $k=2$, we are in the case of Theorem \ref{thm:2gens} and the ideal is in canonical form, and when $k=1$ the ideal has no shared indices and is also in canonical form.
\end{rem}

\begin{proof}
We have
\begin{align*}
    &(x_1y_kg_1,x_2y_1g_2,\ldots,x_ky_{k-1}g_k)\\
    &=(x_1g_1,x_2y_1g_2,\ldots,x_ky_{k-1}g_k) \cap (x_2y_1g_2,\ldots,x_ky_{k-1}g_k,y_k).
\end{align*}
By Lemma \ref{lem:combinecanonical}, the recomposed $\mathfrak{a}$ equals the intersection of the two recomposed pieces, and by Proposition \ref{pr:origgensstay}, we may remove terms divisible by some $x_iy_i$ at any point in the process to get the \fivestep\ of $\mathfrak{a}$.

By Theorem \ref{thm:chain}, the \fivestep\ of the first piece is
\begin{align*}
    (&x_{1}g_1,x_{2}[g_1g_2],x_{3}[g_1g_2g_3],\ldots,x_{{k-1}}[g_1 \cdots g_{k-1}],x_k[g_1 \cdots g_k], \\
    &x_{2}y_{1}g_2,x_{3}y_{1}[g_2g_3],\ldots,x_{{k-1}}y_{1}[g_2\cdots g_{k-1}],x_ky_{1}[g_2\cdots g_k],\\
    &x_{3}y_{2}g_3,\ldots,x_{{k-1}}y_{2}[g_3 \cdots g_{k-1}],x_ky_{2}[g_3 \cdots g_k],\\
    &\vdots \\
    &x_{{k-1}}y_{{k-2}}g_{k-1},x_ky_{{k-2}}[g_{k-1}g_k],\\
    &x_ky_{{k-1}}g_k).
\end{align*}
and the \fivestep\ of the second piece is
\begin{align*}
    (&x_{2}y_{1}g_2,x_{3}y_{1}[g_2g_3],\ldots,x_{{k-1}}y_{1}[g_2\cdots g_{k-1}],y_{1}[g_2\cdots g_k],\\
    &x_{3}y_{2}g_3,\ldots,x_{{k-1}}y_{2}[g_3 \cdots g_{k-1}],y_{2}[g_3 \cdots g_k],\\
    &\vdots \\
    &x_{{k-1}}y_{{k-2}}g_{k-1},y_{{k-2}}[g_{k-1}g_k],\\
    &y_{{k-1}}g_k).
\end{align*}
The \fivestep\ of their intersection is
\begin{align*}
    (&x_1y_kg_1,x_2y_k[g_1g_2],\ldots,x_{k-1}y_k[g_1 \cdots g_{k-1}],\\
    &x_2y_1g_2,x_3y_1[g_2g_3],\ldots,x_ky_1[g_1 \cdots g_k],\\
    &x_3y_2g_3,\ldots,x_ky_2[g_3 \cdots g_k],x_1y_2[g_3 \cdots g_kg_1],\\
    &\vdots \\
    &x_ky_{k-1}g_k,x_1y_{k-1}[g_kg_1],\ldots,x_{k-2}y_{k-1}[g_kg_1 \cdots g_{k-2}]).
\end{align*}
Each term is divisible by a distinct product $x_{i_1}y_{i_2}$, so no term is divisible by another. Hence this is the canonical form of $\mathfrak{a}$. Since $k>2$, this canonical form has more generators than $\mathfrak{a}$.
\end{proof}

\begin{prop}
\label{pr:spreadoutxitoyi}
The \fivestep\ of the polarized neural ideal 
\[\mathfrak{a}=(x_1\cdots x_kg,y_1g_1,\ldots,y_kg_k),\]
where $k>0$ and the only shared indices are the $x_i$ and $y_i$ shown,
is
\begin{align*}
    (&y_1g_1,\ldots,y_kg_k,\\
    &x_1 \cdots x_kg,\\
    &x_2 \cdots x_k[gg_1],x_1x_3 \cdots x_k[gg_2],\ldots,x_1 \cdots x_{k-1}[gg_k],\\
    &x_3 \cdots x_k[gg_1g_2],\ldots,x_1 \cdots x_{k-2}[gg_{k-1}g_k],\\
    &\vdots \\
    &x_k[gg_1 \cdots g_{k-1}],\ldots,x_1[gg_2 \cdots g_k],\\
    &[gg_1 \cdots g_k]).
\end{align*}
Further, this holds with the $x_i$ and $y_i$ reversed.
\end{prop}

\begin{rem}
Unlike in Theorems \ref{thm:chain} and \ref{thm:cycle}, there are many possibilities for the canonical form of $\mathfrak{a}$, depending on which products of $g$ and the $g_j$ have the same LCM's.
\end{rem}

\begin{proof}
We prove this by induction on $k$. If $k=1$, Theorem \ref{thm:2gens} tells us that the \fivestep\ of this ideal is $(x_1g,y_1g_1,[gg_1])$.

Now suppose $k>1$. We have
\begin{align*}
    (x_1\cdots x_kg,y_1g_1,\ldots,y_kg_k) =(x_1,y_1g_1,\ldots,y_kg_k)\cap(x_2\cdots x_kg,y_1g_1,\ldots,y_kg_k).
\end{align*}
By Theorem \ref{thm:2gens} and Lemma \ref{lem:boringgens}, the \fivestep\ of the first piece is $(x_1,g_1,y_2g_2,\ldots,y_kg_k)$. By the induction hypothesis and Lemma \ref{lem:boringgens}, the \fivestep\ of the second piece is
\begin{align*}
    (&y_1g_1,\ldots,y_kg_k, x_2\cdots x_kg, \\
    &x_3 \cdots x_k[gg_2],\ldots,x_2 \cdots x_{k-1}[gg_k],\\
    &x_4 \cdots x_k[gg_2g_3],\ldots,x_2 \cdots x_{k-2}[gg_{k-1}g_k],\\
    &\vdots \\
    &x_k[gg_2 \cdots g_{k-1}],\ldots,x_2[gg_3 \cdots g_k],[gg_1 \cdots g_k]).
\end{align*}
Intersecting these by Lemma \ref{lem:combinecanonical} and removing generators divisible by $x_iy_i$ or that are obviously a multiple of another generator gives the result.

The final statement holds because the $x_i$ and $y_i$ have symmetrical roles in Algorithm \ref{alg:canonform}.
\end{proof}

\begin{prop}
\label{pr:xssplituptoys}
The \fivestep\ of the polarized neural ideal
\[(x_1\cdots x_kg, y_1g_{1},\ldots,y_{\ell}g_{\ell},y_{\ell+1}\cdots y_{t_1}g_{\ell+1},\ldots,y_{t_{s-1}+1} \cdots y_kg_{\ell+s}),\]
where $k \ge 2$, each product involving $g_{\ell+j}$ for $j>0$ has at least 2 $y_i$ as factors, and the only shared indices are the $x_i$ and $y_i$ shown, is
\begin{align*}
    (&y_1g_1,\ldots,y_{\ell}g_{\ell},y_{\ell+1}\cdots y_{t_1}g_{\ell+1},\ldots,y_{t_{s-1}+1}\cdots y_kg_{\ell+s},x_1 \cdots x_kg,\\
    &x_2 \cdots x_k[gg_1],\ldots,x_1 \cdots x_{\ell-1}x_{\ell+1}\cdots x_k[gg_{\ell}],\\
    &x_3\cdots x_k[gg_1g_2],\ldots,x_1 \cdots x_{\ell-2}x_{\ell+1}\cdots x_k[gg_{\ell-1}g_{\ell}],\\
    &\vdots \\
    &x_{\ell} \cdots x_k[gg_1 \cdots g_{\ell-1}],\ldots,x_1x_{\ell+1}\cdots x_k[gg_2 \cdots g_{\ell}],\\
    &x_{\ell+1}\cdots x_k[gg_1 \cdots g_{\ell}]),
\end{align*}
up to removing generators that are a multiple of another generator.
In particular, if $\ell=0$, the ideal is in canonical form already.

Furthermore, this result holds if the roles of the $x_i$ and $y_i$ are reversed.
\end{prop}

\begin{proof}
By Lemma \ref{lem:oneatatime}, it suffices to recompose $\mathfrak{a}$ with respect to the indices $1,\ldots,\ell$. The result now follows from Proposition \ref{pr:spreadoutxitoyi}.
\end{proof}

\begin{example}
The canonical form of $(x_1x_2x_3g_1,y_1y_2g_2,y_3g_3)$ is \[(x_1x_2x_3g_1,y_1y_2g_2,y_3g_3,x_1x_2[g_1g_3]),\]
provided $[g_1g_3] \ne g_1$. If $g_3 \mid g_1$, the canonical form is 
$(y_1y_2g_2,y_3g_3,x_1x_2g_1).$
\end{example}

\begin{example}
We compute the canonical form of the neural ideal 
\[\mathfrak{a}=(x_1x_4x_5,x_2x_3y_1,y_2y_6,y_3y_6,y_3y_4y_5).\]
In order to see how to work with these indices, we picture the ideal as a table, where each column represents a generator and each row a set of indices:
\begin{tabular}{| c | c | c | c | c |}
\hline
$x_1$ & $x_2y_1$ & $y_2$ & - & - \\
\hline
- & $x_3$ & - & $y_3$ & $y_3$ \\
\hline
$x_4x_5$ & - & - & - & $y_4y_5$ \\
\hline
- & - & $y_6$ & $y_6$ & - \\
\hline
\end{tabular}
The shared indices in this ideal are $\{1,2,\ldots,5\}$. However, the indices 4 and 5 are shared by the 1st and 5th generators only, which share 2 indices. So it suffices to consider the indices $1,2,$ and 3.

The ideal $\mathfrak{a}$ has the form $(x_1g_1,x_2y_1g_2,y_2g_3,g_4,g_5)$, where $x_1,x_2,y_1,y_2 \nmid g_j$ for any $j$. Applying Theorem \ref{thm:chain}, we recompose with respect to the indices 1 and 2 to get
\begin{align*}
    (x_1x_4x_5,x_2x_3x_4x_5,x_3x_4x_5y_6, x_2x_3y_1,x_3y_6, y_2y_6, y_3y_6,y_3y_4y_5).
\end{align*}
Removing $x_3x_4x_5y_6$, which is a multiple of $x_3y_6$, we get
\[(x_1x_4x_5,x_2x_3x_4x_5,x_2x_3y_1,x_3y_6,y_2y_6,y_3y_6,y_3y_4y_5).\]

Now the ideal has the form $(g_1,x_3g_2,x_3g_3,x_3g_4,g_5,y_3g_6,y_3g_7)$, where $x_3,y_3 \nmid g_j$ for any $j$. Applying Lemma \ref{lem:lotsofsamexiyi}, we recompose with respect to the index 3 to get:
\begin{align*}
    (x_1x_4x_5,y_2y_6, x_2x_3x_4x_5,x_2x_3y_1,x_3y_6, y_3y_6,x_2x_4x_5y_6,x_2y_1y_6,y_6,\\
    y_3y_4y_5,x_2x_4x_5y_4y_5,x_2y_1y_4y_5,y_4y_5y_6).
\end{align*}
Removing terms divisible by $x_iy_i$ for some $i$ and terms divisible by another term (in particular, divisible by $y_6$), we are left with the canonical form
   \[ (x_1x_4x_5,x_2x_3x_4x_5,x_2x_3y_1,y_6,y_3y_4y_5,x_2y_1y_4y_5).\]
\end{example}

\begin{example}
Let $\mathfrak{a}=(x_1g_1,x_2y_1g_2,x_3y_2g_3,y_3g_4,x_2g_5)$, where the $g_j$ contain no shared indices. This has a chain of $x_i,y_i$ pairs in the first 4 generators, but $x_2$ also appears in the last generator so Theorem \ref{thm:chain} doesn't apply. Theorem \ref{thm:repeats} also doesn't apply since $x_2$ appears with $y_1$ (part of a shared index) in the second generator and alone in the 5th generator. In this case, following Algorithm \ref{alg:shorter} over the set $\{1,2,3\}$
is the right approach to compute this canonical form. We recompose with respect to the index 1 first to get
\[(x_1g_1,x_2y_1g_2,x_3y_2g_3,y_3g_4,x_2g_5,x_2[g_1g_2]).\]
Notice that we now have 3 generators divisible by $x_2$. Recomposing with respect to the index 2, we get
\[(x_1g_1,x_2y_1g_2,x_3y_2g_3,y_3g_4,x_2g_5,x_2[g_1g_2],x_3y_1[g_2g_3],x_3[g_3g_5],x_3[g_1g_2g_3]).\]
Finally, recomposing with respect to the index 3, we get
\begin{align*}
    (&x_1g_1,x_2y_1g_2,x_3y_2g_3,y_3g_4,x_2g_5,x_2[g_1g_2],\\
    &x_3y_1[g_2g_3],x_3[g_3g_5],x_3[g_1g_2g_3],\\
    &y_2[g_3g_4],y_1[g_2g_3g_4],[g_3g_4g_5],[g_1g_2g_3g_4]),
\end{align*}
which is the almost canonical form of the ideal.

To see that this is in fact the right canonical form, one can check each pair of generators $f_1,f_2$ that shares exactly one index $i$ to see that the ideal has another generator dividing $[f_1f_2]/(x_iy_i)$. For example, $x_1g_1$ and $x_2y_1g_2$ share only the index 1, but $x_2[g_1g_2] \mid [x_1g_1 \cdot x_2y_1g_2]/(x_1y_1)$.
\end{example}

\section{Acknowledgments}

We thank Jack Jeffries and Nora Youngs for helpful conversations. In particular, Jack Jeffries inspired Remark \ref{rem:connecttorfstructure}. In addition, Ray Karpman gave extremely helpful combinatorial advice on an earlier version of this paper.

\bibliographystyle{amsalpha}
\bibliography{neuralbib}
\end{document}